\documentclass[12pt]{amsart}
\usepackage{amsfonts, amssymb}
\usepackage{fullpage}
\usepackage{latexsym}

\newcommand\cc{\mathfrak c}
\renewcommand\gg{\mathfrak g}
\newcommand\hh{\mathfrak h}
\newcommand\mm{\mathfrak m}

\newcommand{\inprod}[1]{\langle #1\rangle}

\newcommand\inverse{{^{-1}}}
\renewcommand{\check}{^{\vee}}

\newcommand\ra{\rightarrow}

\DeclareMathOperator{\ad}{ad}
\DeclareMathOperator{\Ad}{Ad}
\DeclareMathOperator{\Aut}{Aut}
\DeclareMathOperator{\Char}{char}

\DeclareMathOperator{\GL}{GL}
\DeclareMathOperator{\Gal}{Gal}
\DeclareMathOperator{\SL}{SL}
\DeclareMathOperator{\PGL}{PGL}

\DeclareMathOperator{\SO}{SO}
\DeclareMathOperator{\SP}{Sp}

\DeclareMathOperator{\Lie}{Lie}

\DeclareMathOperator{\End}{End}

\newcommand{\twobytwo}[4]
 {{\left(\begin{array}{ll} #1 & #2 \\ #3 & #4 \end{array}\right)}}

\numberwithin{equation}{section}

\newtheorem{thm}[equation]{Theorem}

\newtheorem{lem}[equation]{Lemma}
\newtheorem{cor}[equation]{Corollary}
\newtheorem{prop}[equation]{Proposition}

\theoremstyle{definition}
\newtheorem{defn}[equation]{Definition}
\newtheorem{exmp}[equation]{Example}    
\theoremstyle{remark}
\newtheorem{rem}[equation]{Remark}    
\theoremstyle{remark}
\newtheorem{rems}[equation]{Remarks}

\thanks{2000 {\it Mathematics Subject Classification}. 
20G15, 14L24, 
20E42.}

\date{\today}

\title[A Geometric Approach to Complete Reducibility]
{A Geometric Approach to Complete Reducibility}
\author[M.\  Bate]{Michael Bate} 
\address
{School of Mathematics, University of Birmingham, 
Birmingham B15 2TT, UK}
\email{batem@for.mat.bham.ac.uk}
\author[B.\ Martin]{Benjamin Martin}
\address
{Mathematics and Statistics Department,
University of Canterbury, Private Bag 4800,
Christchurch 1, New Zealand}
\email{B.Martin@math.canterbury.ac.nz}
\author[G. R\"ohrle]{Gerhard R\"ohrle}
\address
{School of Mathematics, University of Birmingham, 
Birmingham B15 2TT, UK}
\email{ger@for.mat.bham.ac.uk}

\begin{document}

\begin{abstract}
Let $G$ be a connected reductive linear algebraic group.  
We use geometric methods to investigate $G$-completely 
reducible subgroups of $G$, giving new criteria for $G$-complete reducibility. 
We show that a subgroup of $G$ is $G$-completely reducible if and only if 
it is strongly reductive in $G$; this allows us to use 
ideas of R.W.\ Richardson and Hilbert--Mumford--Kempf 
from geometric invariant theory.
We deduce that a normal 
subgroup of a $G$-completely reducible subgroup of $G$
is again  $G$-completely reducible, thereby providing 
an affirmative answer to a question
posed by J.-P.\  Serre, and conversely we prove that the normalizer of a
$G$-completely reducible subgroup of $G$
is again  $G$-completely reducible.  Some rationality 
questions and applications to the spherical building of $G$ are considered.
Many of our results extend to the case of non-connected $G$.
\end{abstract}

\maketitle

\tableofcontents

\section{Introduction}
\label{s:intro}
Let $G$ be a connected reductive linear algebraic group
defined over an algebraically closed field $k$.  
Following  Serre \cite{serre2}, we say that a (closed) subgroup
$H$ of $G$ is \emph{$G$-completely reducible} ($G$-cr) 
provided that whenever $H$ is contained in a parabolic subgroup $P$ of $G$,
it is contained in a Levi subgroup of $P$; 
for an overview of this concept see for instance \cite{serre1} and 
\cite{serre2}.
In the case $G = \GL(V)$ ($V$ a finite-dimensional $k$-vector space) 
a subgroup $H$ is $G$-cr exactly when 
$V$ is a semisimple $H$-module, 
so this faithfully generalizes
the notion of complete reducibility from representation theory.
The concept of $G$-complete reducibility is part
of the philosophy developed by J.-P.\  Serre, J.\  Tits and others 
to extend standard results from the representation theory of algebraic groups
by replacing representations $H \to \GL(V)$ with 
homomorphisms $H \to G$, 
where the target group is an arbitrary 
reductive algebraic group;
see for instance   \cite{liebeckseitz}, \cite{liebecktesterman}, 
\cite{seitz1}, \cite{serre1.5}, \cite{serre1}, \cite{serre2}, 
and \cite{slodowy}.

In this paper we apply geometric techniques to study $G$-complete 
reducibility.  We are motivated by the philosophy of R.W.~Richardson 
\cite{rich1}.  His insight was that one can study subgroups of $G$ 
indirectly by looking at the action of $G$ on $G^n$ 
by simultaneous conjugation, 
where $n\geq 1$; in this setting, one can apply ideas from 
geometric invariant theory, such as the Hilbert--Mumford Theorem.  
To this end, he introduced the notion of a \emph{strongly reductive} 
subgroup of $G$.  A closed subgroup $H$ of $G$ is said to be  
\emph{strongly reductive in $G$} provided 
$H$ is not contained in any proper parabolic subgroup of 
$C_G(S)$, the centralizer of $S$ in $G$, where
$S$ is a maximal torus of $C_G(H)$, \cite[Def.\  16.1]{rich1}
(this does not depend on the choice of $S$).  He proved that 
if the subgroup $H$ is 
topologically generated by $h_1,\ldots, h_n$, 
then $H$ is strongly reductive in $G$ if and only if the $G$-orbit 
of the $n$-tuple $(h_1,\ldots, h_n)$ 
is closed in $G^n$, \cite[Thm.\  16.4]{rich1}.  
In general, not every $G$-completely reducible 
subgroup is topologically finitely generated, 
but standard arguments (see Remark \ref{rem:topfg} and Lemma \ref{lem:fgapprox}) 
show that one can reduce to this case.

Richardson showed that
a closed subgroup $H$ of $\GL(V)$ 
is strongly reductive if and only if $V$ is a semisimple
$H$-module, \cite[Lem.\  16.2]{rich1}: 
thus strong reductivity and complete reducibility 
are equivalent for subgroups of $\GL(V)$.  Our main result, 
Theorem \ref{mainthm}, asserts that the latter
assertion holds when $\GL(V)$ 
is replaced by an arbitrary reductive group $G$.  This allows 
us to apply results on strong reductivity due to 
Richardson \cite{rich1} and the second author 
\cite{martin1}, \cite{martin2} to study $G$-complete reducibility.  
For example, we deduce immediately that a normal subgroup of 
a $G$-completely reducible subgroup of $G$ is also $G$-completely 
reducible (Theorem \ref{Serrequestion}), thereby answering a 
question of Serre (see  \cite[p.\  24]{serre1}).  Conversely, 
we prove that if a subgroup $H$ is $G$-completely reducible, then 
so is its normalizer $N_G(H)$ (Corollary \ref{cor1converse}).

We then continue our investigations into $G$-complete reducibility, 
centering on the following general question.  Let $f:H\ra G$ be a 
homomorphism of reductive groups and let $K$ be a closed subgroup 
of $H$.  What hypotheses on $H$, $G$, and $f$ guarantee that if 
$K$ is $H$-completely reducible, then $f(K)$ is $G$-completely reducible, 
or vice versa?  Using Lemma \ref{lem:epimorphisms}(ii), 
one can often reduce to the case that $H$ is a closed 
subgroup of $G$.  
Given subgroups $K \subseteq H$ of $G$
with $H$ reductive, 
several of our results in Section \ref{s:cr}
provide criteria that ensure that 
$K$ is $G$-cr if and only if it is $H$-cr
(Corollaries \ref{cor:linearlyreductive} and \ref{cor:levi} 
and Theorem \ref{gcr-regular}).
We also give further conditions to ensure that 
$K$ is $H$-cr, provided it is $G$-cr
(Proposition \ref{ben3}), in particular in case
$(G,H)$ is a reductive pair in the sense of Definition \ref{def:redpair}
(Theorem \ref{red-pair} and Corollary \ref{red-pair-gl}).
When $G=\GL(V)$, the idea is to reduce 
the problem of determining $H$-completely reducible subgroups 
of $H$ to the representation-theoretic problem of determining 
subgroups of $H$ that act completely reducibly on $V$.  
A theorem of Serre (see Theorem \ref{serremain}) gives a 
lower bound on the characteristic $\Char k$ to ensure that 
$H$-complete reducibility is equivalent to $G$-complete 
reducibility for closed subgroups $K$ of $H$.  
We give some related results 
(see Theorem \ref{red-pair}, Corollary \ref{red-pair-gl}, 
Theorem \ref{thm:adjcr},
and the final paragraph of Section \ref{s:cr}), 
which improve on his bound in some circumstances.  
The special case of the adjoint representation is
discussed in greater detail: see Example \ref{ex:adjoint},
Corollary \ref{lie-gcr}, Remarks \ref{rems:lie-gcr}(ii)--(iii), 
and Theorem \ref{thm:adjcr}.
A number of our results give new criteria for 
subgroups to be $G$-completely reducible, such as  
Theorem \ref{Serrequestion}, Theorem \ref{converse}, and 
Proposition \ref{regular-is-cr}.

Much of the previous work on $G$-complete reducibility relies on 
a detailed investigation into the properties of each of the classical 
and exceptional simple algebraic groups (cf.\ \cite{liebeckseitz0}, 
\cite{liebeckseitz}, \cite{liebecktesterman}, and \cite{mcninch}).  
We see our methods as complementing this approach: generally the 
geometric arguments are short and uniform, without requiring a 
case-by-case analysis.

Most of our results can be extended to the case of a non-connected 
group $G$ with $G^0$ reductive; this requires the formalism of 
\emph{R-parabolic subgroups} and \emph{R-Levi subgroups} of $G$, 
to be introduced in Section \ref{s:noncon}.   Often important 
groups associated to a connected group $G$, such as $\Aut G$ 
or centralizers or normalizers of subgroups of $G$, are not connected.  
For example, to prove Theorem 1.1 of \cite{martin1} 
(see Corollary \ref{cor:gcr-finite} below), even just for 
connected groups, one needs to consider non-connected groups.  
Likewise, the formalism of R-parabolic subgroups for non-connected 
groups allows us to deduce Theorem \ref{converse} immediately from 
Propositions \ref{ben3} and \ref{redcent} (see Subsection \ref{sub:cr}).
Nevertheless, to avoid technicalities we have formulated our 
results first for connected $G$: for example, the proof of our main 
result Theorem \ref{mainthm} uses only the standard theory of 
connected reductive groups and their parabolic subgroups.  
The extensions to non-connected $G$, together with the necessary 
preliminary results on R-parabolic subgroups, are postponed to
Section \ref{s:noncon}.  Some of the results in the earlier 
sections, such as Proposition \ref{ben3}, are proved 
using the cocharacter language of Lemma \ref{lem:cochars} to 
facilitate their generalizations in Section \ref{s:noncon}.

\smallskip

The paper is organized as follows.
In Section \ref{s:prelim} 
we first recall some standard
tools from geometric invariant theory we require for 
the sequel.  We give a characterization of parabolic 
subgroups and Levi subgroups of $G$ using cocharacters of $G$.
This is followed by some basic results on 
$G$-complete reducibility and related notions.
We continue by recalling a number of relevant results 
on strongly reductive subgroups.  
This subsection is mostly expository, 
but we believe it is worth reviewing the basics of 
Richardson's theory, because it is not well known.  
The background material is not required for the proof 
of Theorem \ref{mainthm}, with which we start 
Section \ref{s:cr}, but it is needed for the other results that follow.

The heart of the paper is Section  \ref{s:cr} which 
contains most of our results on $G$-cr subgroups,
as outlined above.

In Section \ref{s:building} we discuss  
the building-theoretic approach to
$G$-complete reducibility due to  J.-P.\ Serre, who has shown that 
if $H$ is a $G$-cr subgroup of $G$, then the 
geometric realization of the
fixed point subcomplex of 
the action of $H$ on the Tits building of $G$ is homotopy equivalent to 
a bouquet of spheres. We investigate instances when this
subcomplex is itself a building.

Section \ref{s:rationality} is concerned with rationality 
questions for $G$-cr subgroups. In particular, we 
show that if
$G$ is defined over a perfect field $k$,
then a $k$-subgroup of $G$ is ``$G$-completely reducible over $k$'' 
in an appropriate sense if and only if 
it is $G$-completely reducible.

Finally, in Section \ref{s:noncon} 
we extend some of our results from 
Sections \ref{s:prelim}--\ref{s:rationality} to the case
when $G$ is no longer required to be connected.

Our main references for reductive algebraic groups and their parabolic
subgroups are \cite{borel}, \cite{boreltits}, and \cite{spr2}.

\section{Preliminaries} 
\label{s:prelim}

\subsection{Notation} 
\label{sub:not}
We maintain the notation from the introduction.
In particular, $G$ is a reductive algebraic group defined over 
a field $k$.  Throughout we assume that reductive groups are connected 
(in Section \ref{s:noncon} we use the term non-connected reductive 
group to mean a possibly non-connected group with reductive 
identity component).
With the exception of Section \ref{s:rationality}, 
we assume that $k$ is algebraically closed.
We denote the Lie algebra of $G$ by $\Lie G$ or by $\gg$; 
likewise for closed subgroups of $G$.
For a closed subgroup $H$ of $G$, we denote its 
identity component by $H^0$. 
The centralizer and normalizer of $H$ in $G$ are 
$C_G(H) = \{g\in G \mid ghg\inverse  = h \textrm{ for all } h \in H\}$
and $N_G(H) = \{g\in G \mid gHg\inverse  = H \}$, respectively.
If $S$ is a group acting on $G$ by automorphisms, 
then $C_G(S)$ is the subgroup of $S$-fixed points of $G$.
Also the centralizer of $H$ in $\gg$ is defined by 
$\cc_{\gg}(H) = \{x \in \gg \mid \Ad h(x) = x \textrm{ for all } h \in H\}$.
By a Levi subgroup of $G$ we mean a Levi subgroup of a parabolic 
subgroup of $G$. 
The unipotent
radical of a closed subgroup $H$ of $G$ is denoted by $R_u(H)$.
Let $S$ be a torus of $G$. Then $C_G(S)$ is a Levi subgroup of $G$,
\cite[Thm.\  20.4]{borel}.
Conversely, every Levi subgroup of $G$ is of this form, e.g., see
Lemma \ref{lem:cochars}(ii).
We write $Z(G)$ for the center of $G$.

Let $H$ be a closed (not necessarily connected)  
subgroup of $G$ normalized by some maximal torus $T$
of $G$: that is, 
a \emph{regular} subgroup of $G$ (reductive 
regular subgroups are often also referred to as 
\emph{subsystem subgroups}, e.g.,\ see 
\cite{liebeckseitz0}, \cite{liebeckseitz}).
Let $\Psi = \Psi(G,T)$ denote the set of roots of $G$
with respect to $T$.
In this case the root spaces of $\mathfrak h$ relative to $T$ 
are also root spaces of  $\mathfrak g$ relative to $T$, 
and the set of roots of $H$ with respect to $T$, 
$\Psi(H) = \Psi(H, T) =  \{\alpha \in \Psi \mid \mathfrak g_\alpha
\subseteq \mathfrak h\}$,  is a subset of $\Psi$, 
where $\mathfrak g_\alpha$ denotes the root space in $\mathfrak g$ 
corresponding to $\alpha$.
If $\Char k$ does not divide any of the structure constants
of the Chevalley commutator relations of $G$,
then $\Psi(H)$ is closed under addition in 
$\Psi$ in the sense that if $m,n\in \mathbb N$ 
and $\alpha,\beta\in \Psi(H)$ with $m\alpha+n\beta\in \Psi$, 
then $m\alpha+n\beta\in \Psi(H)$.
If $H$ is reductive and regular, then 
$\Psi(H)$ is a semisimple subsystem of $\Psi$. 

Fix a Borel subgroup $B$ of $G$ containing $T$ and let  
$\Sigma = \Sigma(G, T)$   
be the set of simple roots of $\Psi$ defined by $B$. Then
$\Psi^+ = \Psi(B)$ is the set of positive roots of $G$.
For $\beta \in \Psi^+$ write 
$\beta = \sum_{\alpha \in \Sigma} c_{\alpha\beta} \alpha$
with $c_{\alpha\beta} \in \mathbb N_0$.
A prime $p$ is said to be \emph{good} for $G$
if it does not divide $c_{\alpha\beta}$ for any $\alpha$ and $\beta$, 
\cite{SS}.

Recall that a linear algebraic group $S$, not necessarily connected, 
is said to be {\em linearly reductive} if every rational representation 
of $S$ is semisimple.  It is well known that in characteristic zero, 
$S$ is linearly reductive if and only if $S^0$ is reductive.  
In characteristic $p>0$, $S$ is linearly reductive if and only 
if every element of $S$ is semisimple if and only if $S^0$ 
is a torus and $|S/S^0|$ is coprime to $p$, see \cite[\S4, Thm.\  2]{nagata}.

\subsection{Characteristic zero}
\label{sub:charzero}

The notions of strong reductivity in $G$ and
$G$-complete reducibility are uninteresting in characteristic zero, 
as a closed subgroup $H$ is strongly reductive in $G$ if and only if 
$H^0$ is reductive if and only if $H$ is $G$-cr (cf.\ Lemma \ref{linred-cr}).
We therefore assume for the remainder of the paper
that $k$ has characteristic $p >0$.
Nevertheless, all of our results hold in characteristic zero
with the obvious modifications.

\subsection{Geometric invariant theory}
\label{sub:git}

We recall some results from geometric invariant theory 
required in the sequel, see \cite[Ch.\ 3]{newstead}, \cite[\S 2]{BaRi}.
Let $G$ be a reductive group acting on an affine variety $X$ (we assume 
all actions are left actions).

\begin{defn}
\label{defn:stability}
For $x \in X$ let $G \cdot x$ denote the $G$-orbit of $x$ in $X$ and 
$C_G(x)$ the stabilizer of $x$ in $G$.  
Let $Z = \bigcap_{x \in X}C_G(x)$ be the kernel of the action
of $G$ on $X$.
Following \cite[1.4]{rich1} we say that $x \in X$ is a
\emph{stable point} for the action of $G$ or a 
\emph{$G$-stable point} provided the orbit $G \cdot x$ is closed in $X$ and
$C_G(x)/Z$ is finite. 
\end{defn}

\begin{defn}
\label{defn:limit}
Let $\phi : k^* \to X$ be a morphism of algebraic varieties. 
We say that  $\underset{t\to 0}{\lim}\, \phi(t)$ exists 
if there exists a morphism $\widehat\phi :k \to X$ 
(necessarily unique) whose restriction to $k^*$ 
is $\phi$; if this limit exists, then we set
$\underset{t\to 0}{\lim}\, \phi(t) = \widehat\phi(0)$.
\end{defn}

By $Y(G)$ we denote the set of all cocharacters 
$\lambda : k^* \to G$ of $G$.  If $\lambda\in Y(G)$ 
and $g\in G$, then we define $g\cdot \lambda \in Y(G)$ by 
$(g\cdot \lambda)(t) = g\lambda(t)g^{-1}$; 
this gives a left action of $G$ on $Y(G)$.
It follows easily from Definition \ref{defn:limit}
that if $\underset{t\to 0}{\lim}\,\lambda(t)\cdot x$ 
exists for a cocharacter $\lambda \in Y(G)$, 
then this limit belongs to the closure $\overline{G\cdot x}$ 
of $G\cdot x$ in $X$.  
The following result, known as the
Hilbert--Mumford Theorem \cite[Thm.\ 1.4]{kempf}, gives a converse to this.

\begin{thm} 
\label{thm:HMT}
Let $G$ be a reductive group acting on an affine variety $X$, 
and let $x\in X$.  For any $y$ in the closure of $G\cdot x$, there 
exists $\lambda \in Y(G)$ such that 
$\underset{t\to 0}{\lim}\,\lambda(t) \cdot x$ exists and
belongs to $G \cdot y$.
\end{thm}

An important tool in the geometric approach to $G$-complete reducibility 
is a strengthened version of the Hilbert--Mumford 
Theorem due to Kempf \cite[Thm.\ 3.4]{kempf}: 
roughly, this says that if $y$ belongs to the complement
$\overline{G\cdot x} \setminus G\cdot x$,
then there is a \emph{canonical} way of choosing a cocharacter $\lambda$ such 
that $\underset{t\to 0}{\lim}\,\lambda(t) \cdot x$ lies in $G \cdot y$,
a so-called ``optimal'' $\lambda$.  
We refer to this theorem and its corollaries as the 
Hilbert--Mumford--Kempf Theorem.  
It is the main ingredient in the proof of 
Proposition \ref{martinthm1} and underlies the rationality result
\cite[Thm.\ 4.2]{kempf}, which is used indirectly in our proof of 
Theorem \ref{thm:perfectcr}.

We require the characterization of parabolic subgroups of $G$
in terms of cocharacters of $G$, 
see \cite[2.1--2.3]{rich1} and 
\cite[Prop.\  8.4.5]{spr2}:

\begin{lem} 
\label{lem:cochars} 
Given a parabolic subgroup $P$ of $G$ and any Levi subgroup $L$ of $P$,
there exists $\lambda \in Y(G)$ 
such that the following hold:
\begin{itemize}
\item[(i)]   $P = P_\lambda := \{g\in G \mid \underset{t\to 0}{\lim}\,
        \lambda(t) g \lambda(t)^{-1} \textrm{ exists}\}$.
\item[(ii)]  $L = L_\lambda := C_G(\lambda(k^*))$.
\item[(iii)] The map $c_\lambda : P_\lambda \to L_\lambda$ defined by
\[
c_\lambda(g) := \underset{t\to 0}{\lim}\, \lambda(t)g \lambda(t)\inverse
\]
        is a surjective homomorphism of algebraic groups.
        Moreover, $L_\lambda$ is the set of fixed points of $c_\lambda$
        and $R_u(P_\lambda)$ is the kernel of $c_\lambda$.
\end{itemize}
Conversely, given  any $\lambda \in Y(G)$ 
the subset $P_\lambda$ defined as in part (i) 
is a parabolic subgroup of $G$,
$L_\lambda$ is a Levi subgroup of $P_\lambda$, and
the map $c_\lambda$ as defined in part (iii) 
has the described properties.  
Moreover, $P_\lambda$ is a proper subgroup if and only if 
$\lambda(k^*)\not\subseteq Z(G)$.
\end{lem}

Let $H$ be a reductive subgroup of $G$.
There is a natural inclusion $Y(H) \subseteq Y(G)$ of
cocharacter groups.
If necessary, 
we distinguish between the ambient groups we are working
in by writing $P_\lambda(H)$, $P_\lambda(G)$, etc., 
for $\lambda \in Y(H)$; usually we just write
$P_\lambda$ rather than $P_\lambda(G)$.
It is obvious from Lemma~\ref{lem:cochars}
that if $\lambda \in Y(H)$, then $P_\lambda(H) = P_\lambda(G)\cap H$ 
and similarly
for $L_\lambda(H)$ and $R_u(P_\lambda(H))$.
We record this in our next result.

\begin{cor}
\label{cor:parsofsubgp}
Let $H$ be a reductive subgroup of $G$.
If $Q$ is a parabolic
subgroup of $H$ and $M$ is a Levi subgroup of $Q$, 
then there exists a parabolic subgroup $P$
of $G$ and a Levi subgroup $L$ of $P$ such that 
$Q = P\cap H$, $M = L\cap H$ and $R_u(Q) = R_u(P) \cap H$.
\end{cor}

We say that a closed subgroup $H$ of $G$ is \emph{topologically 
generated by $h_1,\ldots, h_n \in G$} provided
that $H$ is the Zariski closure of the subgroup of $G$ generated by 
these elements. 
As we are concerned with topologically finitely generated subgroups of $G$,
we are interested in the diagonal action of $G$ on the affine variety $G^n$ 
(for some $n\in {\mathbb N}$) by simultaneous conjugation:
\[
g\cdot (g_1,\ldots, g_n) := (gg_1g^{-1},\ldots, gg_ng^{-1}). 
\]
Let $(h_1,\ldots, h_n)\in G^n$ and let $H$ be the subgroup of $G$ 
topologically generated by $h_1,\ldots, h_n$.  
Let $\lambda \in Y(G)$.  
It follows easily from Definition \ref{defn:limit} 
that $\underset{t\to 0}{\lim}\,\lambda(t) \cdot (h_1,\ldots, h_n)$ 
exists if and only if $\underset{t\to 0}{\lim}\,\lambda(t)h_i\lambda(t)^{-1}$ 
exists for $1\leq i\leq n$, and this is the case if and only if 
all of the $h_i$ lie in $P_\lambda$, cf.~Lemma~\ref{lem:cochars}(i).  
Since the $h_i$ 
topologically generate $H$, this occurs if and only if 
$H$ is a subgroup of $P_\lambda$, and in this case, 
$c_\lambda(h_1),\ldots, c_\lambda(h_n)$ 
topologically generate $c_\lambda(H)$, where 
$c_\lambda$ is the map defined in 
Lemma \ref{lem:cochars}(iii).

\subsection{Basic properties of $G$-cr, $G$-ir, and $G$-ind subgroups} 
\label{sub:basics-gcr}
The concept of $G$-complete reducibility is 
a relative notion
depending on the embedding of the subgroup into $G$.  Note that 
$G$ is trivially a $G$-cr subgroup of itself.
If $H$ is a closed $G$-cr subgroup of $G$, 
then $H^0$ is reductive,  \cite[Property 4]{serre1}.
In characteristic zero the converse holds.
This follows from a well-known 
result due to G.\  Mostow \cite{Mostow}.
More generally, we have the following result 
(e.g., see \cite[Lem.\  11.24]{Jantzen}):  

\begin{lem}
\label{linred-cr}
Let $S$ be a linearly reductive subgroup of $G$.
Then $S$ is $G$-completely reducible.
\end{lem}

An important class of $G$-cr subgroups consists of those that are 
not contained in any proper parabolic subgroup of $G$ at all 
(they are trivially $G$-cr). 
Following Serre, we call them
\emph{$G$-irreducible} ($G$-ir), \cite{serre2}.  
As with the concept of $G$-complete reducibility, 
this terminology stems from the fact that 
it coincides with the usual notion of irreducibility
in the classical case $G = \GL(V)$.
Observe that every overgroup of a $G$-irreducible subgroup of $G$ 
is itself $G$-irreducible. 
Moreover, following Serre \cite{serre2}, we say that a subgroup $H$ of $G$ is 
\emph{$G$-indecomposable} ($G$-ind) provided that $H$ is not contained in 
any Levi subgroup of any proper parabolic 
subgroup of $G$. 
Again, in the classical case $G = \GL(V)$ this coincides with 
the usual property  of $V$ being an indecomposable $H$-module.

Note that a closed subgroup $H$ of $G$ is strongly 
reductive in $G$ if and only if $H$ is $C_G(S)$-ir, 
where $S$ is a maximal torus of $C_G(H)$.

Part (i) of our next result is due to Slodowy 
in the special case $G = \GL(V)$, \cite[Lem.\ 11]{slodowy2}.

\begin{cor}
\label{cor:ben2}
Let $K \subseteq H$ be closed subgroups of $G$ with 
$H$ reductive.
\begin{itemize}
\item[(i)] If $K$ is $G$-irreducible, then it is $H$-irreducible.
\item[(ii)] If $K$ is $G$-indecomposable, then it is $H$-indecomposable.
\end{itemize}
\end{cor}

\begin{proof}
It follows directly from Corollary \ref{cor:parsofsubgp}
that if $K$ is contained in a proper parabolic subgroup of $H$, 
then it is contained in a proper parabolic subgroup of $G$ as well,
giving (i).
Similarly, if $K$ is contained in a proper Levi subgroup of $H$, 
then it is also contained in a proper Levi subgroup of $G$,
again by Corollary \ref{cor:parsofsubgp}, 
giving (ii). 
\end{proof}

\begin{rems}
(i).
If $H$ is not $G$-irreducible, then taking $K = H$ gives a trivial 
counterexample to the converse of Corollary \ref{cor:ben2}(i).
Example \ref{ex:serrecountereg} below 
(taking $K= \PGL_p(k)$, $H=\PGL_m(k)$, and $G=\GL(\hh)$) 
shows that the converse of 
Corollary \ref{cor:ben2}(i) can fail even when $H$ is $G$-ir.
Likewise, if $H$ is not $G$-ind, taking $K = H$ again 
gives a trivial counterexample to the converse of 
Corollary \ref{cor:ben2}(ii).

(ii).
Serre's Theorem \ref{serremain}(ii) below
shows that if $\Char k$ is sufficiently large and $G = \GL(V)$, then 
the implication of Corollary \ref{cor:ben2}(i) does hold 
for $G$-cr and $H$-cr in place of $G$-ir and $H$-ir, respectively.
However, this statement does not hold in general, e.g., see 
Example \ref{ex:liebeckcountereg} below.
Several of the results in Section \ref{s:cr}
provide criteria to ensure the analogous implication holds 
for $G$-complete reducibility and $H$-complete reducibility, for instance see
Proposition \ref{ben3}, Corollary \ref{cor:linearlyreductive}, 
Theorem \ref{gcr-regular}, Theorem \ref{red-pair}, and 
Corollary \ref{red-pair-gl}.
\end{rems}

\begin{rem}
\label{rem:topfg}
When considering questions of $G$-complete reducibility, etc., 
it is technically convenient to work with topologically 
finitely generated subgroups of $G$.  Unfortunately, not every 
reductive subgroup of $G$ has this property: for example, if $k$ 
is the algebraic closure of the prime field $\mathbb {F}_p$,
then every finitely generated subgroup of $G$ is finite.
The following argument of Richardson \cite[Prop.\  16.9]{rich1}
allows us to reduce to the case of
topologically finitely generated subgroups of $G$
in our study of $G$-cr subgroups; it is fundamental in 
employing the methods from geometric invariant theory 
outlined in Subsection \ref{sub:git}.
By Theorem \ref{thm:algclosrat} and 
Remark \ref{rem:algclosratirind} below, one may assume without loss that 
$k$ is  transcendental over $\mathbb F_p$.
Then one can show that any reductive subgroup $H$ of $G$ is
topologically finitely generated, cf.\ \cite[Lem.\ 9.2]{martin1}.  
We use this idea repeatedly in what follows.
\end{rem}

In the proof of Theorem \ref{thm:perfectcr} we 
need another method for reducing to the topologically finitely 
generated case.

\begin{lem}
\label{lem:fgapprox}
Let $H$ be a closed subgroup of $G$.
There exists a finitely 
generated subgroup $\Gamma$ of $H$ such 
that for any parabolic subgroup $P$ 
of $G$ and any Levi subgroup $L$ of 
$P$, $P$ contains $H$ if and only if 
$P$ contains $\Gamma$, and $L$ contains 
$H$ if and only if $L$ contains $\Gamma$.  
\end{lem}

\begin{proof}
It is well known that $G$ has only finitely many conjugacy classes of 
parabolic subgroups, and each parabolic subgroup
$P$ has exactly one $P$-conjugacy class of Levi 
subgroups, so we can choose a finite set of representatives 
$P_1,\ldots, P_m$ and $L_1,\ldots, L_n$ for the set of $G$-conjugacy 
classes of parabolic subgroups and Levi subgroups respectively.  
For any $i = 1,\ldots, m$ and $j = 1,\ldots, n$ and any $H' \subseteq H$, 
set $C_i(H')=\{g \in G \,|\, H' \subseteq gP_ig^{-1}\}$ and 
$D_j(H')=\{g \in G \,|\, H' \subseteq gL_jg^{-1}\}$.  
Each $C_i(H')$ and $D_j(H')$ is closed, and if 
$H' \subseteq H''$, then $C_i(H'') \subseteq C_i(H')$ 
and $D_j(H'') \subseteq D_j(H')$.  
The descending chain condition on 
closed subsets of $G$, together with a simple application of 
Zorn's Lemma, implies that for some 
finitely generated subgroup $\Gamma \subseteq H$, 
we have $C_i(\Gamma)=C_i(H)$ and $D_j(\Gamma)=D_j(H)$
for every $i$ and $j$.  The result now follows.
\end{proof}

We want to investigate how the  properties of
$G$-complete reducibility, $G$-irreducibility, 
and $G$-indecomposability behave under 
homomorphisms of the ambient groups, cf.\ \cite[Cor.\ 4.3]{serre2}.
The next result and part (ii) of Lemma \ref{lem:epimorphisms}
answer this question for epimorphisms.

\begin{lem}
\label{lem:isogenynew}
 Let $f:G_1\ra G_2$ be an isogeny of reductive groups and let 
$\lambda\in Y(G_1)$.  Set $\mu=f\circ \lambda \in Y(G_2)$.  Then
 \begin{itemize}
  \item[(i)] $f(P_\lambda)=P_\mu$, $f(L_\lambda)=L_\mu$;
  \item[(ii)] $f^{-1}(P_\mu)=P_\lambda$, $f^{-1}(L_\mu)=L_\lambda$;
  \item[(iii)] $L_\lambda$ is a proper subgroup of $G_1$ if and only if $L_\mu$ is a proper subgroup of $G_2$;
  \item[(iv)] $f(R_u(P_\lambda))= R_u(P_\mu)$.
 \end{itemize}
\end{lem}

\begin{proof}
 Suppose that we have a connected 
subset $S$ of $G_1$ and an element $x\in G_1$ such that $f(x)$ 
centralizes $f(S)$.  Then $[x,S]\subseteq \ker{f}$, and we deduce 
that $[x,S] = \{ 1 \}$, as $\ker{f}$ is finite.  
This implies that 
$f(L_\lambda)=f(C_{G_1}(\lambda(k^*))) = C_{G_2}(f(\lambda(k^*))) = 
C_{G_2}(\mu(k^*)) = L_\mu$, and part (iii) also follows. 
 
If $B$ is a Borel subgroup of $G_1$, then 
$f(B)$ is a Borel subgroup of $G_2$ (\cite[Prop.\ 16.7]{borel}).  
This implies that $f(P_\lambda)$ is a parabolic subgroup of $G_2$.
It is clear from the definition of a limit, Definition \ref{defn:limit},  
that $f(P_\lambda)\subseteq P_\mu$.  As $f(P_\lambda)$ contains 
$f(L_\lambda)=L_\mu$, a Levi subgroup of $P_\mu$, 
\cite[Prop.\ 4.4(c)]{boreltits} implies that 
$f(P_\lambda)=P_\mu$, which completes the proof of part (i).

Next we observe  that $\ker f$ is central in $G_1$, so it is contained 
in any Levi subgroup of $G_1$.  Thus (ii) follows easily from (i).

It is clear that $f(R_u(P_\lambda))$ is a closed connected unipotent 
normal subgroup of $P_\mu$, so $f(R_u(P_\lambda))\subseteq R_u(P_\mu)$.  
Since $L_\lambda$ and $L_\mu$ are Levi subgroups of $P_\lambda$ and $P_\mu$, 
respectively, part (i) together with a simple dimension-counting argument 
implies that $f(R_u(P_\lambda))= R_u(P_\mu)$, as required.
\end{proof}

Note that if $f:G_1\ra G_2$ is an 
isogeny of reductive groups, then 
any $\mu\in Y(G_2)$ is of the form 
$f\circ \lambda$ for some $\lambda\in Y(G_1)$. 

Let $f:G_1\rightarrow G_2$ be a homomorphism of algebraic groups.  
We say that $f$ is \emph{non-degenerate} provided 
$({\rm ker}f)^0$ is a torus, cf.\ \cite[Cor.\ 4.3]{serre2}.
In particular, $f$ is non-degenerate if $f$ is an isogeny.

\begin{lem}
\label{lem:epimorphisms}
Let $G_1$ and $G_2$ be reductive groups.
\begin{itemize}
\item[(i)]
Let $H$ be a 
closed subgroup of $G_1\times G_2$.  
Let $\pi_i:G_1\times G_2\rightarrow G_i$ 
be the canonical projection for $i = 1,2$.  Then $H$ 
is $(G_1\times G_2)$-cr 
($(G_1\times G_2)$-ir, $(G_1\times G_2)$-ind) 
if and only if $\pi_i(H)$ is $G_i$-cr ($G_i$-ir, $G_i$-ind) for $i=1,2$.
\item[(ii)]
Let $f:G_1\rightarrow G_2$ be an epimorphism.  
Let $H_1$ and $H_2$ be closed subgroups of $G_1$ and $G_2$, respectively.
\begin{itemize}
\item[(a)] 
If $H_1$ is $G_1$-cr ($G_1$-ir, $G_1$-ind), 
then $f(H_1)$ is $G_2$-cr ($G_2$-ir, $G_2$-ind).
\item[(b)] 
If $f$ is non-degenerate, then $H_1$ is $G_1$-cr ($G_1$-ir, $G_1$-ind) 
if and only if $f(H_1)$ is $G_2$-cr ($G_2$-ir, $G_2$-ind), 
and $H_2$ is $G_2$-cr ($G_2$-ir, $G_2$-ind) if and only if 
$f^{-1}(H_2)$ is $G_1$-cr ($G_1$-ir, $G_1$-ind).
\end{itemize}
\end{itemize}
\end{lem}

\begin{proof}
(i). The parabolic subgroups of $G_1\times G_2$ 
are precisely the subgroups of the form $P_1\times P_2$, 
where $P_i$ is a parabolic subgroup of $G_i$.  
If $P_1\times P_2$ is such a subgroup, then the 
Levi subgroups of $P_1\times P_2$ are precisely 
the subgroups of the form $L_1\times L_2$, where 
$L_i$ is a Levi subgroup of $P_i$.  Part (i) now follows.

(ii). Let $N=\ker f$.
It is standard that there exists a closed reductive subgroup $M_1$ of 
$C_{G_1}(N^0)$ with  $M_1 \cap N^0$ finite and $G_1=M_1 N^0$ 
(cf.\ the proof of Lemma \ref{lem:normalcomplement}).  
It follows 
from Lemma \ref{lem:isogenynew}(i)--(iii) 
that the required result holds for isogenies, so we can assume 
that $N$ is connected,
$G_1 = M_1\times N$, and $f$ is the projection from 
$G_1$ to $M_1$.  Parts (a) and (b) now follow from part (i); 
note that if $N$ is 
a torus, then any closed subgroup of $N$ is trivially $N$-cr, 
$N$-ir, and $N$-ind.
\end{proof}

We require the following characterization of 
topologically finitely generated $G$-irreducible subgroups 
of $G$ in terms of stability due to R.W.\ Richardson, 
\cite[Prop.\  16.7]{rich1}. 

\begin{prop}
\label{rich-stable}
Let $x_1, \ldots, x_n \in G$ 
and  let $H$ be the 
subgroup of $G$ that is topologically generated by  $x_1, \ldots, x_n$. 
Then $H$ is $G$-irreducible if and only if 
$(x_1, \ldots, x_n)$ is a stable point under the diagonal action of $G$
on $G^n$, 
that is if and only if the
$G$-orbit of $(x_1, \ldots, x_n)$ in $G^n$ is closed and 
$C_G(H)^0 = Z(G)^0$. 
\end{prop}

Observe that the last equivalence 
in Proposition \ref{rich-stable} follows from the definition of stability, 
Definition \ref{defn:stability},  as $C_G(H)$ is the stabilizer in $G$ of 
$(x_1, \ldots, x_n) \in G^n$.

\subsection{Strongly reductive subgroups of $G$}
\label{sub:strong}

If $H$ is a strongly reductive subgroup 
of $G$, then $H^0$ is reductive, 
\cite[Lem.\  16.3]{rich1}, or \cite[\S 6]{martin1}.
In characteristic zero, the converse also holds, \cite[\S 16]{rich1}, or
\cite[Prop.\ 6.6]{martin1}.
However, in positive characteristic, this is a more 
subtle notion, which depends on the embedding of $H$ into $G$.
Trivially, $G$ is a strongly reductive subgroup of itself.
We require several results on strong reductivity, 
the first of which 
is due to the second author, \cite[Thm.\ 2]{martin2}.

\begin{prop}
\label{martinthm1}
Let $H$ be a strongly reductive subgroup of $G$ and let $N$ be a 
closed normal subgroup of $H$. Then $N$ is strongly reductive in $G$.
\end{prop}

The following two results are due to 
R.W.\  Richardson \cite[Prop.\  16.9, Thm.\  16.4]{rich1}.

\begin{prop}
\label{Rich-linearlyreductive}
Let $S$ be a linearly reductive group acting on the reductive group $G$ 
by automorphisms and let $H= C_G(S)^0$.
Suppose $K$ is a closed subgroup of $H$.
Then $K$ is strongly reductive in $H$ if and only if
it is strongly reductive in $G$.
\end{prop}

Recall that $H = C_G(S)^0$ is reductive, \cite[Prop.\  10.1.5]{rich0}.

The following result 
gives a geometric interpretation for topologically finitely generated
algebraic subgroups of $G$ that are strongly reductive in $G$.

\begin{prop}
\label{Rich-closedorbit}
Let $x_1, \ldots, x_n \in G$ 
and  let $H$ be the 
subgroup of $G$ that is topologically generated by  $x_1, \ldots, x_n$. 
Then $H$ is strongly reductive in $G$ if and only if 
the $G$-orbit of $(x_1, \ldots, x_n)$ under the diagonal action of $G$
on $G^n$ is closed.
\end{prop}

Observe that the case $n=1$ is simply the characterization of 
semisimple elements in $G$, \cite[Cor.\ 3.6]{SS}.

Using the Hilbert--Mumford Theorem \ref{thm:HMT}  and
the map $c_\lambda$ from Lemma \ref{lem:cochars}(iii), 
we give another characterization of strong reductivity.

\begin{lem}
\label{lem:srcrit}
Let $H$ be a closed subgroup of $G$.  Then $H$ is 
strongly reductive in $G$ if and only if for every 
cocharacter $\lambda$ of $G$ with  $H\subseteq P_\lambda$, 
there exists $g\in G$ 
such that $c_\lambda(h)=ghg^{-1}$ for every $h\in H$.
\end{lem}

\begin{proof}
For simplicity, we assume that $H$ is topologically finitely generated, 
say by $h_1,\ldots, h_n$.  (This suffices for the application to 
Theorem \ref{thm:adjcr} below; the general case follows from a 
slightly more complicated argument involving Lemma \ref{lem:fgapprox}.)  
For any cocharacter $\lambda$ such 
that $H\subseteq P_\lambda$, the subgroup $c_\lambda(H)$ is topologically 
generated by the elements of the tuple 
$(h_1',\ldots, h_n'):=\underset{t\to 0}{\lim}\,\lambda(t)
\cdot (h_1,\ldots, h_n)$, 
and this tuple belongs to the closure of the $G$-orbit
$G\cdot (h_1,\ldots, h_n)$ in $G^n$.  

If $H$ is strongly reductive in $G$, then $G\cdot (h_1,\ldots, h_n)$ 
is closed, by Proposition \ref{Rich-closedorbit}, so 
$(h_1',\ldots, h_n')=g\cdot (h_1,\ldots, h_n)$ for some $g\in G$, 
so $c_\lambda(h_i)=h_i'=gh_ig^{-1}$ for every $i$.  
This implies that $c_\lambda(h)=ghg^{-1}$ for every $h\in H$.  

Conversely, suppose that for every $\lambda \in Y(G)$ with 
$H\subseteq P_\lambda$, there exists $g\in G$ with
$c_\lambda(h)=ghg^{-1}$ for every $h\in H$.  Then for every 
such cocharacter $\lambda$, 
there exists $g\in G$ such that 
$\underset{t\to 0}{\lim}\, 
\lambda(t)\cdot (h_1,\ldots, h_n)=g\cdot (h_1,\ldots, h_n)$.  
The Hilbert--Mumford Theorem \ref{thm:HMT} 
now implies that $G\cdot (h_1,\ldots, h_n)$ 
is closed, so $H$ is strongly reductive in $G$, again by  
Proposition \ref{Rich-closedorbit}.
\end{proof}

\begin{rem}
The conclusion of Lemma \ref{lem:srcrit}, that $H$ and 
$c_\lambda(H) \subseteq L_\lambda$ are
$G$-conjugate, does not imply \emph{a priori} 
that $H$ is itself in a Levi subgroup of
$P_\lambda$. 
That this is indeed the case is 
the content of the reverse implication of 
Theorem \ref{mainthm} below.
\end{rem}

\subsection{Non-closed subgroups}
\label{sub:closing}

Observe that the notions of strong reductivity in $G$, 
$G$-complete reducibility, etc., apply also to non-closed subgroups.
Clearly, a subgroup $H$ of $G$ is $G$-cr if and only if 
its Zariski closure $\overline H$ is, and similarly for  
$G$-irreducibility, $G$-indecomposability and strong reductivity in $G$.
For convenience we only consider closed subgroups of $G$ in the sequel; 
note that all the results we give hold for non-closed subgroups as well, 
wherever this makes sense.

\section{$G$-Complete Reducibility}
\label{s:cr}

\subsection{$G$-complete reducibility and strong reductivity in $G$}
\label{sub:gcr-sr}

The basis for our study of $G$-completely reducible groups is 

\begin{thm} 
\label{mainthm}
Let $G$ be reductive and suppose $H$ is a closed subgroup of $G$.
Then $H$ is $G$-completely reducible if and only if 
$H$ is strongly reductive in $G$.
\end{thm}

\begin{proof}
Suppose that $H$ is $G$-cr, and let $S$ be a maximal torus of $C_G(H)$.  
Now suppose that $H$ is contained in some proper parabolic 
subgroup  $Q$ of $C_G(S)$.  
There exists a parabolic subgroup 
$P$ of $G$ such that $Q=C_G(S)\cap P$ (\cite[Prop.\  4.4(c)]{boreltits}).  
Note that $S\subseteq P$, because $S$ is central in $C_G(S)$.  
Since $H$ is $G$-cr, there is a Levi subgroup  $L$ of $P$ such that 
$H\subseteq L$, and $T:=Z(L)^0$ is a torus of $C_P(H)$.  
Now $S$ is a maximal torus of $C_P(H)$, so $gTg^{-1}\subseteq S$ 
for some $g\in C_P(H)$.  It follows that 
$C_G(S)\subseteq C_G(gTg^{-1})= gC_G(T)g^{-1}=gLg^{-1}\subseteq P$, 
whence $Q=C_G(S)$, a contradiction.  Thus $H$ is strongly reductive 
in $G$.

If $H$ is strongly reductive in $G$ and 
$S$ is a maximal torus of $C_G(H)$, 
then $H$ is not in any proper parabolic subgroup of $C_G(S)$.
Since $S$ is a torus, $L:=C_G(S)$ is a Levi subgroup of $G$.
Let $Q$ be a parabolic subgroup of $G$ containing $L$ as a Levi subgroup.
Then, since $H$ is in no proper parabolic subgroup of $L$,
it follows from \cite[Prop.\  4.4(c)]{boreltits} 
that $Q$ is minimal among all parabolic subgroups containing $H$.
Now let $P$ be a parabolic subgroup of $G$
containing $H$. Then $H \subseteq P \cap Q$.  
If $P'$ is a parabolic subgroup of $G$ with 
$P'\subseteq P$ and $M'$ is a Levi subgroup 
of $P'$, then there exists a Levi subgroup $M$ 
of $P$ such that $M'\subseteq M$.
Therefore, we may assume 
that $P$ is minimal 
subject to $P \supseteq H$.
Since $(P \cap Q) R_u(Q)$ is a parabolic subgroup of $G$ 
(\cite[Prop.\  4.4(b)]{boreltits})
contained in $Q$, the minimality of $Q$ implies that  
$Q = (P \cap Q) R_u(Q)$.
By \cite[Prop.\  4.4(b)]{boreltits} 
it follows that $P$ contains a Levi subgroup, $M_Q$ say, 
of $Q$.
By symmetry, 
$Q$ contains a Levi subgroup, $M_P$ say, of $P$.
It follows that $P \cap Q$ contains a common Levi subgroup of both 
$P$ and $Q$. For, 
fix Levi subgroups, $L_P$ and $L_Q$, of $P$ and $Q$ respectively such that
$L_P\cap L_Q$ contains a maximal torus of $G$. 
Then we have a decomposition
\begin{equation}
\label{fac}
P\cap Q =(L_P \cap L_Q)(L_P \cap R_u(Q))(R_u(P) \cap L_Q)(R_u(P) \cap R_u(Q))
\end{equation}
(this is standard; see Lemma \ref{lem:nonconn}(iii) below for a proof).  
Moreover, $R_u(P \cap Q)$ is the product of the last three factors
in \eqref{fac}.
Since $M_P$ is reductive, $M_P \cap R_u(P \cap Q)$ is trivial.
Thus there is a bijective homomorphism from $M_P$ onto 
a subgroup of $L_P \cap L_Q$.
Since $L_P$ and $M_P$ are $P$-conjugate, we get
$L_P \subseteq L_Q$. Likewise, $L_Q \subseteq L_P$. 
It follows that $M := L_P = L_Q$ is a Levi subgroup of both
$P$ and $Q$, as claimed.
 
Let $P^-$ be the unique parabolic subgroup of $G$ opposite to 
$P$, such that $P \cap P^- = M$.
We may factor $R_u(Q)$ as follows:
\begin{equation}
\label{dec}
R_u(Q) = (R_u(Q) \cap R_u(P^-))(R_u(Q) \cap R_u(P)).
\end{equation}
For the set of roots of $G$ with respect to some maximal torus 
of $G$ in $M$ decomposes
as a disjoint union $\Psi(R_u(P)) \cup \Psi(M) \cup \Psi(R_u(P^-))$.
The decomposition \eqref{dec} follows, as $R_u(Q) \cap M$ is trivial.

Now, since $L$ and $M$ are Levi subgroups of $Q$, 
there exists $x$ in $R_u(Q)$ with $xMx\inverse = L$. 
Thanks to \eqref{dec}, we may write $x=yz$
with $y \in R_u(Q) \cap R_u(P^-)$ and $z \in R_u(Q) \cap R_u(P)$.
Then, since $z \in P \cap Q$, we have $zMz\inverse \subseteq P \cap Q$. 
Thus we may assume  that $z=1$.
Now, as $y \in P^-$, we see that $L = yMy\inverse$ lies in $P^-$. 
Thus $H$ lies in $P^-$.
Consequently, $H \subseteq P \cap P^- = M$.
It follows that $H$ is $G$-cr, as required.
\end{proof}

\begin{rem}
\label{rem:nonreduced}
Serre has observed that the definitions of 
strong reductivity in $G$ and $G$-complete reducibility make 
sense for an arbitrary, possibly non-smooth, subgroup subscheme 
of $G$.  The proof of Theorem \ref{mainthm} goes through virtually 
unchanged in this setting; note that if $P$ and $Q$ are parabolic 
subgroups of $G$, then the scheme-theoretic intersection of $P$ 
and $Q$ is smooth, see \cite[vol.\ 3, Lem.\ 4.1.1]{SGA}.
\end{rem}

\begin{cor}
\label{mainthm-cor}
Let $H$ be a closed subgroup of $G$.
Then the following are equivalent:
\begin{itemize}
\item[(i)] $H$ is strongly reductive in $G$;
\item[(ii)] $H$ is $G$-completely reducible;
\item[(iii)] $H$ is $C_G(S)$-irreducible, where $S$ is a maximal torus of 
$C_G(H)$;
\item[(iv)] for every parabolic subgroup $P$ of $G$
which is minimal with respect to
containing $H$, the subgroup 
$H$ is $L$-irreducible for some Levi subgroup $L$ of $P$;
\item[(v)] there exists a 
parabolic subgroup $P$ of $G$
which is minimal with respect to
containing $H$, such that 
$H$ is $L$-irreducible for some Levi subgroup $L$ of $P$.
\end{itemize}
\end{cor}

\begin{proof}
The equivalences between (i), (ii), and (iii) follow from 
Theorem \ref{mainthm} and the definition of strong reductivity.  
It is clear that (iv) implies (v).  
Given a maximal torus $S$ of $C_G(H)$, $C_G(S)$ is a Levi subgroup of
some parabolic subgroup $P$ of $G$; if $H$ is $C_G(S)$-ir, 
then $P$ is minimal with respect to containing $H$, 
by \cite[Prop.\ 4.4(c)]{boreltits}, so (iii) implies (v).  
If there exists a parabolic subgroup $P$ minimal with respect to 
containing $H$ and $H$ is an $L$-ir subgroup of some Levi subgroup 
$L$ of $P$, then the proof of Theorem \ref{mainthm} implies that $H$ 
is $G$-cr, so (v) implies (ii).  Finally, if $H$ is $G$-cr, then $H$ 
is contained in a Levi subgroup $L$ of $P$, where $P$ is any parabolic 
subgroup minimal with respect to containing $H$.  
By \cite[Prop.\ 4.4(c)]{boreltits}, $H$ is $L$-ir, so (ii) implies (iv).  
\end{proof}

\begin{rem}
\label{reduce-to-ir}
By Corollary \ref{mainthm-cor} 
the study of $G$-cr subgroups of $G$ reduces to the
study of $L$-ir subgroups of the Levi subgroups $L$ of $G$
(including the case $L=G$).
\end{rem}

Using Theorem \ref{mainthm} and the results on strong reductivity
from Section \ref{s:prelim},
we immediately deduce 
results on $G$-complete reducibility.
Our next result, which follows directly from 
Proposition \ref{Rich-closedorbit} 
and Theorem \ref{mainthm},
allows us to use methods from geometric
invariant theory to study $G$-completely reducible subgroups.
It is crucial for a number of results to follow.

\begin{cor}
\label{cor:closedorbit}
Let $x_1, \ldots, x_n \in G$ (for some $n \in \mathbb N$) and 
let $H$ be the
subgroup of $G$ topologically generated by $x_1, \ldots, x_n$. 
Then $H$ is $G$-completely reducible if and only if 
the orbit of $(x_1, \ldots, x_n)$ under the diagonal action of $G$
on $G^n$ is closed.
\end{cor}

Observe that the notions of $G$-complete reducibility, etc.,  
all apply to finite subgroups of $G$.
Our next result follows immediately from Theorem  \ref{mainthm} and
\cite[Thm.~1.2]{martin1}, which is the corresponding result for 
strongly reductive subgroups of $G$.

\begin{cor}
\label{cor:gcr-finite}
There is only a finite number $c_N$ of $G$-conjugacy classes of 
$G$-completely reducible subgroups of fixed order $N$.
\end{cor}

\begin{rems}
\label{rem:conjclass}
(i). In \cite[Prop.\  2.1]{LMS} it is proved, for $G$ 
simple and of adjoint type, that there is a uniform bound on $c_N$  
as in Corollary \ref{cor:gcr-finite} that depends only on $N$ and the type of 
$G$, but not on the field $k$.

(ii). Observe that in general there exist infinitely many
$G$-conjugacy classes of  
connected $G$-cr subgroups of bounded dimension.  For example, 
the subtori $T_n$ of $\GL_2(k)$ defined by
$T_n := 
\left.\left\{ \twobytwo{t^n}{0}{0}{t}\right| t\in k^* \right\}$ 
for $n$ a positive integer
are pairwise non-conjugate.  For a more 
subtle example involving semisimple subgroups, 
see \cite[Cor.\  4.5]{liebecktesterman}.
\end{rems}

\subsection{Normal subgroups of $G$-cr subgroups}
\label{sub:Serrequestion}

Our next result, which follows immediately from Proposition \ref{martinthm1}
and Theorem \ref{mainthm}, gives an affirmative answer to 
Serre's question, \cite[p.\  24]{serre1}.
The special case when $G = \GL(V)$ is just a particular 
instance of Clifford Theory.

\begin{thm} 
\label{Serrequestion}
Let $H$ be a closed subgroup of $G$ with closed normal subgroup $N$.
If $H$ is $G$-completely reducible, then so is $N$.
In particular, if $H$ is $G$-completely reducible, then so is $H^0$.
\end{thm}

\begin{rem}
\label{Serre-converse}
Serre proves a converse to Theorem \ref{Serrequestion}
in \cite[Property 5]{serre1} under the assumption that 
the index of $N$ in $H$ is prime to $\Char k = p$.
Examples show that this restriction 
cannot be removed. 
For instance, let $U$ be a non-trivial finite unipotent subgroup of $G$.
Then, by a construction due to Borel and Tits
\cite[\S 3]{boreltits2}, there
exists a parabolic subgroup $P$ of $G$ such that 
$U \subseteq R_u(P)$. In particular, $U$ is not $G$-cr, but clearly
$U^0 =\{1\}$ is.
In Theorem \ref{converse} below we give a converse of 
Theorem \ref{Serrequestion} 
without characteristic restrictions but with the  
additional assumption that 
$H$ contains the connected centralizer of $N$ in $G$.
\end{rem}

\subsection{Normalizers and centralizers of $G$-cr subgroups}
\label{sub:geometry}

If $\Char k = 0$ and $H$ is a reductive subgroup of $G$, then
$C_G(H)^0$ and $N_G(H)^0$ are also known to be reductive. 
In positive characteristic this is still true provided 
$H$ is linearly reductive, thanks to \cite[Prop.\ 10.1.5]{rich0}.
However,  it is clearly false in general.
For instance, if $V$ is an indecomposable, non-simple 
$H$-module with isomorphic top and socle, then 
the centralizer of $H$ in $\GL(V)$ is not reductive.  
Nevertheless, as an application of Corollary \ref{cor:closedorbit}, 
we have the following result.

\begin{prop}
\label{redcent}
Let $H$ be a 
$G$-completely reducible
subgroup of $G$. Then $C_G(H)^0$ is  reductive.
Moreover, let $K$ be a closed subgroup of $G$ satisfying
$H^0C_G(H)^0 \subseteq K \subseteq N_G(H)$. Then $K^0$ is reductive.
In particular, $N_G(H)^0$ is reductive.
\end{prop}

\begin{proof}
By Remark \ref{rem:topfg}, we may assume that
$H$ is topologically finitely generated, say by 
$x_1, \ldots, x_n \in G$.
Observe that the stabilizer of $(x_1, \ldots, x_n)$ in 
$G$ is just $C_G(H)$.
As $H$ is $G$-completely reducible, 
the $G$-orbit of $(x_1, \ldots, x_n)$ 
in $G^n$ is closed,
by Corollary \ref{cor:closedorbit}, 
whence $G/C_G(H)$ is affine (e.g., \cite[Lem.\ 10.1.3]{rich0}).
Thanks to 
\cite[Thm.\ A]{rich}, it follows that $C_G(H)^0$ is  reductive.

Now $H^0 \times C_G(H)^0$ is reductive and thus so is
its image $H^0C_G(H)^0$ in $G$.
Since $N_G(H)$ is a finite extension of $H^0C_G(H)^0$,
see  \cite[Lem.\  6.8]{martin1},
the desired result follows.
\end{proof}

\begin{rem}
\label{red-cent-ls}
The main result in \cite[Thm.\  1]{liebeckseitz0} 
asserts that if $p > 7$ 
and $H$ is any reductive subgroup of the simple exceptional 
group $G$, then $H$ is $G$-cr.
In fact, the bounds established in loc.\ cit.\ are
much more detailed and smaller depending on the types of 
$G$ and $H$. Using this result, Liebeck and Seitz
proceed  in \cite[Thm.\  2]{liebeckseitz0} to show 
that the connected centralizer $C_G(H)^0$ 
for each such $H$ is again reductive.
Proposition \ref{redcent}  gives a short proof of
this fact for any reductive $G$. 
\end{rem}

Theorem \ref{converse} 
can be viewed as a  partial converse of 
Theorem \ref{Serrequestion} 
complementing Serre's result \cite[Property 5]{serre1},
cf.\  Remark \ref{Serre-converse}.

\begin{thm}
\label{converse}
Let $H$ be a $G$-completely reducible
subgroup of $G$ and suppose $K$ is a closed subgroup of $G$ satisfying
$HC_G(H)^0 \subseteq K \subseteq N_G(H)$.
Then $K$ is $G$-completely reducible.
\end{thm}

\begin{proof}
Let $P$ be a parabolic subgroup of $G$ containing $K$.
Since $H$ is $G$-cr, there is a Levi subgroup $L$ of $P$ with $L \supseteq H$.
Thanks to Lemma \ref{lem:cochars} there exists 
a cocharacter $\lambda$ of $G$ such that
$P = P_\lambda$ and $L = L_\lambda$.
Since $\lambda(k^*) \subseteq C_G(H)^0 \subseteq K$, we see that
$K$ is $c_\lambda$-stable, where
$c_\lambda$ is the homomorphism defined in Lemma \ref{lem:cochars}(iii). 
We have
\begin{equation}
\label{K}
K = (K \cap L)(K\cap R_u(P)).
\end{equation}
This follows, as every $x \in K \subseteq P = LR_u(P)$ has a unique
factorization $x = x_1x_2$ with $x_1 \in L$  and $x_2 \in R_u(P)$.
Since $K$ is $c_\lambda$-stable, we have 
$c_\lambda(x) = x_1 \in K\cap L$ and this implies \eqref{K}.
For any $x\in K\cap R_u(P)$ we
have a morphism $\phi_x:k \to K\cap R_u(P)$ given by
$\phi_x(t) = \lambda(t) x \lambda(t)\inverse$ for $t\in k^*$
and $\phi_x(0) = c_\lambda(x) = 1$.
Thus the image of $\phi_x$ is a connected subvariety of $K\cap R_u(P)$ 
containing $1$ and $x$.  We deduce that $K\cap R_u(P)$ is connected.
By assumption $HC_G(H)^0 \subseteq K$ and thus $H^0C_G(H)^0 \subseteq K$.
Thus Proposition \ref{redcent} implies that $K\cap R_u(P)$ is trivial.
Consequently, thanks to \eqref{K}, we have $K \subseteq L$, as required.
\end{proof}

The following are immediate consequences of 
Theorems \ref{Serrequestion}  and \ref{converse}.

\begin{cor}
\label{cor1converse}
Let $H$ be a closed subgroup of $G$.
Then $H$ is $G$-completely reducible if and only if 
$N_G(H)$ is. 
\end{cor}

\begin{cor}
\label{cor2converse}
Let $H$ be a closed subgroup of $G$.
If $H$ is $G$-completely reducible, then so is $C_G(H)$.
\end{cor}

\begin{proof}
Since $C_G(H)$ is normal in $N_G(H)$, the result follows from
Theorem \ref{Serrequestion} and Corollary \ref{cor1converse}.
\end{proof}

We observe that in the  classical  case $G = \GL(V)$,
Corollaries \ref{cor1converse} and \ref{cor2converse}
are just consequences of Clifford Theory and Wedderburn's Theorem.

In general, the converse of Corollary \ref{cor2converse}
is false: e.g., let $H$ be a Borel subgroup of $G$. 
However, we do have the following partial converse.  
The proof follows immediately from Lemma \ref{linred-cr} and 
Lemma \ref{lem:irrcentlinred} below.  

\begin{cor}
\label{kempf1}
Let $H$ be a closed subgroup of $G$.
If $C_G(H)$ is $G$-irreducible, then 
$H$ is linearly reductive.  
In particular, $H$ is $G$-completely reducible.
\end{cor}

The hypothesis of Corollary \ref{kempf1} is very restrictive:  
since $H \subseteq C_G(C_G(H))$,  
Remark \ref{rem:topfg} and 
Proposition \ref{rich-stable} imply that 
$H^0 \subseteq Z(G)^0$.

Our next result is similar in nature to Theorem \ref{converse};
it also provides
another criterion to ensure that if $K \subseteq H \subseteq G$
with $H$ reductive and $K$ is $G$-cr, then $K$ is also 
$H$-cr.

\begin{prop}
\label{ben3}
Let $K\subseteq H$ be closed subgroups of $G$ with $H$ reductive.
Suppose that $H$ contains a maximal torus of $C_G(K)$ and that  
$K$ is $G$-completely reducible.
Then $K$ is $H$-completely reducible and $H$ is $G$-completely reducible.
\end{prop}

\begin{proof}
Let $S$ be a maximal torus of $C_G(K)$ contained in $H$.
Then $S$ is also a maximal torus of $C_H(K)$.
Since $K$ is $G$-cr, $K$ is $C_G(S)$-ir, 
by Corollary \ref{mainthm-cor}.
Applying Corollary \ref{cor:ben2}(i) to 
$K\subseteq C_H(S) \subseteq C_G(S)$, we see that
$K$ is $C_H(S)$-ir.
Thus $K$ is $H$-cr, again by Corollary \ref{mainthm-cor}.

Let $P$ be a parabolic subgroup of $G$ containing $H$.
Then $P$ also contains $K$ and $S$.
Since $K$ is $G$-cr and $S\subseteq P$, it follows that $C_G(S)$ 
is contained in a Levi subgroup of $P$ (cf.\ the first paragraph 
of the proof of Theorem \ref{mainthm}): say $C_G(S) \subseteq L_\lambda$, 
where $\lambda\in Y(G)$ and $P=P_\lambda$.
We have 
$\lambda(k^*) \subseteq C_G(C_G(S))^0 = Z(C_G(S))^0$.
Clearly, $S \subseteq Z(C_G(S))^0$. Since 
$K \subseteq C_G(S)$, we have $ Z(C_G(S))^0 \subseteq C_G(K)$.
Because $S$ is a maximal torus of $C_G(K)$, we
get the equality $S = Z(C_G(S))^0$.
So we have  $\lambda(k^*)  \subseteq S$ which is contained in $H$, 
by assumption; so  $\lambda$ is a cocharacter of $H$.
But now we see that $P_\lambda(H) = P_\lambda \cap H = H$.
It follows from Lemma \ref{lem:cochars} that  $\lambda(k^*) \subseteq Z(H)$.
So $H  \subseteq  C_G(\lambda(k^*))  = L_\lambda$, as desired.
\end{proof}

\subsection{$G$-complete reducibility and regular subgroups}
\label{sub:regular}

Recall from Subsection \ref{sub:not} that $H$ is a regular subgroup of 
$G$ provided it is normalized by a maximal torus of $G$.

\begin{prop}
\label{regular-is-cr}
Let $H$ be a regular reductive subgroup of $G$. Then $H$ is 
$G$-completely reducible.
\end{prop}

\begin{proof}
Let $T$ be a maximal torus of $G$
normalizing $H$.  If $T \subseteq H$, then,
since $C_G(T) = T$ and $T$ is $G$-cr by Lemma \ref{linred-cr},
Proposition \ref{ben3} implies 
that $H$ is $G$-cr. 
In the general case $HT$ is $G$-cr by the argument just given; 
thus $H$ is $G$-cr by Theorem \ref{Serrequestion}.
\end{proof}

Our next result follows 
immediately from Proposition \ref{Rich-linearlyreductive}
and Theorem \ref{mainthm}.

\begin{cor}
\label{cor:linearlyreductive}
Let $S$ be a linearly reductive group acting on the reductive group $G$ 
by automorphisms and let $H= C_G(S)^0$.
Suppose $K$ is a closed subgroup of $H$.
Then $K$ is $H$-completely reducible if and only if
it is $G$-completely reducible.
\end{cor}

Note that $H = C_G(S)^0$ is reductive, thanks to 
\cite[Prop.\  10.1.5]{rich0}: in fact, 
by Corollary \ref{cor:linearlyreductive}, $H$ is $G$-cr, because 
$H$ is $H$-cr.
In the examples below we indicate how Corollary \ref{cor:linearlyreductive} 
leads to new criteria for completely reducible subgroups.

The following special case of Corollary \ref{cor:linearlyreductive},
when $S$ is a subtorus of $G$ so that 
$H = C_G(S)$ is a Levi subgroup of $G$, is  
a result due to J.-P.\  Serre \cite[Prop.\  3.2]{serre2}. 
This is again part 
of the philosophy mentioned in the Introduction, because 
for $G = \GL(V)$, the statement simply reduces to the fact that 
$V$ is a semisimple $K$-module if and only if 
every $K$-submodule of $V$ 
is semisimple.

\begin{cor}
\label{cor:levi}
Let $K$ be a closed subgroup of a Levi subgroup $L$ of $G$. 
Then $K$ is $L$-completely reducible if and only if $K$ is 
$G$-completely reducible.
\end{cor}

Other typical applications of Corollary \ref{cor:linearlyreductive}
are when $S$ is the group generated by a graph automorphism of $G$, 
or $S$ is the group generated by a semisimple element of $G$ such that
$C_G(S)^0$ is a subgroup of $G$ of maximal semisimple rank.
The subsystems corresponding to maximal semisimple rank
subgroups  of a simple group $G$ are determined by means of the algorithm of 
Borel and de Siebenthal \cite{BoSe}, see also 
\cite[Ex.\ Ch.\ VI \S 4.4]{bou}.
We give some examples.

\begin{exmp}
\label{gcr-classical}
Suppose that $p \ne 2$.
Let $V$ be a finite-dimensional $k$-vector space.
Let $H$ be either $\SP(V)$ or $\SO(V)$ and let 
$K$ be a closed subgroup of $H$. 
Observe that $H$ is the fixed point subgroup of 
an involution of $\GL(V)$ and thus Corollary \ref{cor:linearlyreductive} applies.
Then $K$ is $H$-completely reducible if and only if
it is $\GL(V)$-completely reducible, that is, 
if and only if $V$ is a semisimple $K$-module;
see also \cite[Cor.\  16.10]{rich1}, and 
\cite[Ex.\  3.2.2(b)]{serre2}.
\end{exmp}

\begin{exmp}
\label{g2}
Suppose that $p > 3$.
Let $G = \SO_8(k)$ and let $H$ be the subgroup of $G$ of type $G_2$  
that is the fixed point subgroup
of $G$ under the triality graph automorphism of $G$.
Let $K$ be a closed subgroup of $H$.
Then $K$ is $H$-completely reducible if and only if
it is $G$-completely reducible.
By our assumption on $p$ the natural 
$8$-dimensional module $V_G(8)$ of $G$ 
splits as an $H$-module into a direct
sum $V_G(8)|_H = V_H(7) \oplus k$
of the simple $7$-dimensional $H$-module $V_H(7)$ and the trivial module.
It follows from 
Corollary \ref{cor:linearlyreductive} and 
Example \ref{gcr-classical} that
$K$ is $H$-completely reducible if and only if
$V_H(7)$ is a semisimple $K$-module.
We observe that this equivalence does also hold in characteristic $3$, 
see \cite[Ex.\  3.2.2(c)]{serre2}. 
\end{exmp}

\begin{exmp}
\label{e8}
Let $G$ be of type $E_8$
and suppose that $p \ne 5$.
Let $H$ be a maximal rank subgroup of $G$ of type $A_4A_4$.
Let $K$ be a closed subgroup of $H$.
Then $K$ is $H$-completely reducible if and only if
it is $G$-completely reducible.
This follows 
from Corollary \ref{cor:linearlyreductive}, as 
$H$ is the centralizer  in $G$ of an element of order $5$.
Using Lemma \ref{lem:epimorphisms}
this can then be interpreted in terms of the semisimplicity of 
the natural modules of the $A_4$-factors of $H$ for the corresponding 
projections of $K$ into these factors. 
\end{exmp}

More generally, we have the following result.

\begin{thm}
\label{gcr-regular}
Suppose that $p$ is good for $G$.
Let $H$ be a regular reductive subgroup of $G$. Let $K$ be 
a closed subgroup of $H$. Then 
$K$ is $H$-completely reducible if and only if 
$K$ is $G$-completely reducible.
\end{thm}

\begin{proof}
Let $T$ be a maximal torus that normalizes $H$.  
There exists 
a subtorus $S$ of $T$ such that $S$ centralizes 
$H$, $H\cap S$ is finite and $T\subseteq HS$.  
Applying Lemma \ref{lem:epimorphisms}(i) and (ii) 
to the product map $H\times S\rightarrow HS$, 
we see that a closed subgroup $K$ of $H$ is $H$-cr if and 
only if it is $HS$-cr.  Thus we may assume that $H$ contains $T$.

As $H$ is a regular reductive subgroup of $G$, 
it is $G$-cr, by Proposition \ref{regular-is-cr}.
Thanks to
Corollaries \ref{mainthm-cor} and \ref{cor:levi}, 
we may assume that $H$ is $G$-ir.
By Remark \ref{rem:topfg},  
we may assume that $H$ is topologically finitely generated.
Then, since $C_G(H)^0 = Z(G)^0$ by Proposition \ref{rich-stable}, 
we see that $H$ has maximal semisimple rank.
Finally, the result follows from the algorithm of 
Borel and de Siebenthal (\cite{BoSe} or \cite[Ex.\ Ch.\ VI \S 4.4]{bou}),
Carter's criterion \cite[Prop.\ 11]{carter},
and  Corollary \ref{cor:linearlyreductive}.
\end{proof}

\subsection{$G$-complete reducibility, separability and reductive pairs}
\label{sub:redpair}

Now we consider the 
interaction of subgroups of $G$ with the Lie algebra $\Lie G = \gg$ of $G$.

\begin{defn}
\label{def:separable} 
We say a closed subgroup $H$ of $G$ is \emph{separable in $G$}
if the Lie algebra centralizer $\cc_\gg(H)$ of $H$ 
equals the Lie algebra of $C_G(H)$ (that is, if the scheme-theoretic centralizer of $H$ in $G$ is smooth).  If the former properly contains 
the latter, then we say that $H$ is \emph{non-separable in $G$}.
\end{defn}

\begin{exmp}
\label{exmp:allsep}
Any closed subgroup $H$ of $G=\GL(V)$ is separable in $G$.  
For separability means precisely that the centralizers 
of $H$ in $\GL(V)$ and in $\Lie \GL(V)$ 
have the same dimension, and for $\GL(V)$ this holds, because
the centralizer of $H$ in $\GL(V)$ is the open 
subset of invertible elements of the centralizer of $H$
in $\Lie \GL(V) \cong \End V$.
\end{exmp}

\begin{exmp}
\label{exmp:PGL2}
(Cf.\ \cite[Ex.\ 3.4(b)]{martin3}).  
Let $p=2$,  $G = \PGL_2(k)$ and let $T$ be a maximal torus of $G$.  
Then $N_G(T)$ is non-separable in $G$.
\end{exmp}

\begin{rem}
\label{rem:seprem}
 If $G$ itself is non-separable in $G$, 
there exists a reductive group $\widehat{G}$ 
such that $\widehat{G}$ is generated by $G$ and a central torus, 
and $\widehat{G}$ is separable in $\widehat{G}$, 
see \cite[Thm.\ 4.5]{martin3}.
For example, $G = \SL_p(k)$ is non-separable in $G$, 
and we can take $\widehat{G} = \GL_p(k)$.
\end{rem}

\begin{rem}
\label{rem:separable} 
The terminology in Definition \ref{def:separable} 
is motivated as follows. Suppose that $H$ is topologically generated
by $x_1,\ldots ,x_n$ in $G$.
Then the orbit map  $G \ra G \cdot (x_1,\ldots ,x_n)$ 
is separable if and only if 
\[
\cc_\gg(H) = \cc_\gg(\{x_1,\ldots ,x_n\}) = 
\Lie C_G((x_1,\ldots ,x_n)) = \Lie C_G(H)
\] 
(cf. \cite[Prop.\  6.7]{borel}), i.e., if and only if $H$ is separable in $G$.
\end{rem}

\begin{defn}
\label{def:redpair} 
Following Richardson \cite{rich2}, we call 
$(G,H)$ a \emph{reductive pair} provided $H$ is a closed reductive subgroup of 
$G$ and $\Lie G$ decomposes as an $H$-module into a direct sum
\[
\Lie G = \Lie H \oplus \mm, 
\]
where $H$ acts via the adjoint action $\Ad_G$.
\end{defn}

For a list of examples of reductive pairs we refer 
to P.\ Slodowy's article  \cite[I.3]{slodowy}.  
The next example gives a further class of 
reductive pairs.

\begin{exmp}
\label{exmp:redpair}
Let $H$ be a reductive subgroup of $G$ containing a maximal
torus $T$ of $G$, so that in particular, $H$ is regular reductive. 
One can show that if 
$p$ does not divide any of the structure constants
of the Chevalley commutator relations of $G$,
then $(G,H)$ is a reductive pair.
The complement to $\Lie H$ in $\Lie G$ is the 
sum of the root spaces of $\Lie G$ corresponding to the roots in 
$\Psi(G,T)$ that lie outside $\Psi(H,T)$.

This is no longer valid if we relax the restriction on $p$.
For instance, suppose $p=2$ and
let $G$ be of type $B_2$. If $H$ is the regular reductive subgroup 
of type $A_1^2$ generated by the short root subgroups of $G$, 
then  $\Lie H$ does not admit an $H$-stable complement in $\Lie G$.
\end{exmp}

The following observation, due to Serre, 
gives many more examples of reductive pairs. 

\begin{rem}
\label{rem:redpair}
Let  $f : G_1 \to G_2$ be a homomorphism of reductive groups
and let $df : \Lie G_1 \to  \Lie G_2$ be the induced homomorphism on 
the Lie algebras.
Suppose that there exists a symmetric Ad-invariant bilinear form 
(``Ad-invariant form'' for short) 
$(\cdot \,, \cdot)_2$ on 
$\Lie G_2$ which is non-degenerate.
We then define an Ad-invariant form on $\Lie G_1$ via
$(x, y)_f := (df(x), df(y))_2$ for $x,y \in \Lie G_1$.
If $(\cdot \,, \cdot)_f$ 
is non-degenerate and $df:\Lie G_1 \to \Lie f(G_1)$ is surjective, 
then $(G_2, f(G_1))$ is a reductive pair. 
To see this, 
we take $\mm$ to be the  orthogonal complement 
of $\Lie f(G_1)$ in $\Lie G_2$ with respect 
to $(\cdot \,, \cdot)_2$, an $f(G_1)$-stable 
subspace of $\Lie G_2$.  Our hypotheses imply 
that the restriction of $(\cdot \,, \cdot)_2$ 
to $\Lie f(G_1)$ is non-degenerate, so $\Lie f(G_1)\cap {\mm}=\{0\}$.  
As $(\cdot \,, \cdot)_2$ is non-degenerate, we have 
$\dim{\mm}+\dim{f(G_1)}= \dim{G_2}$, as required.

Suppose that $\Lie G_1$ is simple and admits 
a non-degenerate Ad-invariant form $(\cdot \,, \cdot)_1$.  
It can be shown that $(\cdot \,, \cdot)_1$ is the unique 
Ad-invariant form up to scalar multiplication, so we have 
$(\cdot \,, \cdot)_f = \delta_f (\cdot \,, \cdot)_1$ for some 
$\delta_f \in k$.  Thus $(G_2,f(G_1))$ is a reductive pair 
as long as $\delta_f \neq 0$.
In analogy to \cite[\S 2]{dynkin} we 
call $\delta_f$ the \emph{Dynkin index} of $f$.  
For tables of the Dynkin index of the fundamental representations 
see \cite[\S 5]{MPR}, where this invariant is 
called the ``second index''; here $(\cdot \,, \cdot)_1$ 
and $(\cdot \,, \cdot)_2$ are fixed by appropriate 
normalization conditions.
\end{rem}

\begin{thm}
\label{red-pair}
Suppose that $(G,H)$ is a reductive pair.  
Let $K$ be a closed subgroup of $H$ 
such that $K$ is a separable subgroup of $G$.  
If $K$ is $G$-completely reducible, 
then it is also $H$-completely reducible.
\end{thm}

\begin{proof}
By Remark \ref{rem:topfg},
we can assume that $K$ is topologically 
finitely generated, say by $k_1,\ldots ,k_n$.
Let $C$ be the $G$-orbit of  $(k_1,\ldots, k_n)$ in $G^n$.
By assumption, the orbit map $G \ra C$ 
is separable, cf.~Remark \ref{rem:separable}.  
Thanks to a generalization of a standard 
tangent space argument of Richardson \cite[Thm.\ 4.1]{rich2} 
to this situation, see 
Slodowy \cite[Thm.\  1]{slodowy}, the intersection 
$C \cap H^n$ 
is a finite union of $H$-conjugacy classes, each of which is
closed in $C \cap H^n$;  
the second of these assertions
follows directly from the proof in loc.\ cit., as 
each irreducible component of $C \cap H^n$ is a single $H$-orbit.

Now suppose that $K$ is $G$-cr.  Then
$C$ is closed in $G^n$ by Corollary \ref{cor:closedorbit}, and so
the $H$-orbit of $(k_1,\ldots ,k_n)$ 
is closed in $H^n$ by the above argument.  
Thus, using Corollary \ref{cor:closedorbit} again, we see that $K$ 
is $H$-cr, as desired.
\end{proof}

Example \ref{exmp:allsep} shows that the separability 
hypothesis is automatically satisfied for the case 
$G=\GL(V)$.  We obtain an immediate consequence of 
Theorem \ref{red-pair}, which 
is in the spirit of Serre's Theorem \ref{serremain}(ii) below.

\begin{cor}
\label{red-pair-gl}
Suppose that $(\GL(V),H)$ is a reductive pair
and  $K$ is a closed subgroup of  $H$. 
If  $V$ is a semisimple $K$-module, then
$K$ is $H$-completely reducible.
\end{cor}

In our next example we look at the special case of the adjoint representation. 

\begin{exmp}
\label{ex:adjoint}
Let $H$ be a simple group of adjoint type and let 
$G = \GL(\Lie H)$.
We have a symmetric non-degenerate Ad-invariant bilinear form 
on $\Lie G \cong \End (\Lie H)$ 
given by the usual trace form 
and its restriction to $\Lie H$ is just the Killing form of $\Lie H$.
Since $H$ is adjoint and $\Ad$ is a closed embedding, 
$\ad : \Lie H \to \Lie \Ad(H)$ is surjective. 
Thus it follows from the arguments in the first paragraph of 
Remark \ref{rem:redpair} that if the Killing form of $\Lie H$ 
is non-degenerate, then $(G,H)$ is a reductive pair.

Suppose first that $H$ is a simple classical group of adjoint type
and $p > 2$.
The Killing form 
is non-degenerate for ${\mathfrak {sl}}(V)$, ${\mathfrak {so}}(V)$, or
${\mathfrak {sp}}(V)$  if and only if $p$  does not
divide  $2\dim V$, $\dim V-2$, or $\dim V+2$, respectively,
cf.\ \cite[Ex.~Ch.~VIII, \S 13.12]{bou2}.  
In particular, for  $H$ adjoint of type 
$A_n$, $B_n$, $C_n$, or $D_n$, the Killing form 
is non-degenerate if $p > 2$  and $p$ does not
divide $n+1$, $2n-1$, $n+1$, or $n-1$, respectively. 

Now suppose that $H$ is a simple exceptional group of adjoint type. 
If $p$ is good for $H$, then
the Killing form of $\Lie H$ is non-degenerate; this 
was noted by Richardson (see \cite[\S 5]{rich2}).
Thus if $p$ satisfies the appropriate condition, then 
$(\GL(\Lie H),H)$ is a reductive pair and Corollary \ref{red-pair-gl}
applies.
\end{exmp}

\subsection{$G$-ir subgroups}
\label{sub:gir}
We can say more about the centralizers of $G$-irreducible subgroups
than in Corollary \ref{cor2converse}:

\begin{lem}
\label{lem:irrcentlinred}
Let $H$ be a $G$-irreducible subgroup of $G$.  
Then $C_G(H)$ is linearly reductive.
\end{lem}

\begin{proof}
Suppose that $u$ is a non-trivial unipotent 
element of $C_G(H)$.  Then $H$ centralizes the 
non-trivial unipotent subgroup generated by $u$.
By a construction due to Borel and Tits
\cite[Prop.\  5.4(b)]{martin1}, $H$ is not $G$-irreducible, 
a contradiction.  
We conclude that all elements of $C_G(H)$ are semisimple, 
and we deduce from \cite[\S4, Thm.\  2]{nagata} that $C_G(H)$ 
is linearly reductive.
\end{proof}

In \cite[Thm.\ 3]{liebeckseitz0} M.~Liebeck and G.~Seitz 
proved that if $G$ is simple of exceptional type, $H$ is a 
simple subgroup of $G$, and $p >7$, then $H$ 
is a separable subgroup of $G$.  The next result 
shows that even in low characteristic there is 
not much freedom for a $G$-ir subgroup to be non-separable in $G$.

\begin{prop}
\label{prop:irrnonsep}
Let $H\subseteq G$ be a non-separable 
$G$-irreducible subgroup of $G$.  Then there 
exists a regular non-separable 
subgroup $M$ of $G$ containing $H$ with $M^0$ reductive.
\end{prop}

\begin{proof}
Let $x\in \cc_\gg(H)$ such that $x$ 
is not in the Lie algebra of $C_G(H)$, and set $M=C_G(x)$.
Since $H \subseteq M$ and $H$ is $G$-ir, $M$ is $G$-ir;  
in particular, $M^0$ is reductive.  Thanks to \cite[Lem.~2.4]{martin3},
$x$ is semisimple, so $M$ contains a maximal torus of $G$, 
which implies that $M$ is regular.  
Since $C_G(M)\subseteq C_G(H)$, we see that
$x$ is not in the Lie algebra of $C_G(M)$, but 
by definition of $M$ we have $x \in \cc_\gg(M)$, 
so $M$ is non-separable in $G$.
\end{proof}

The following result answers a question raised by 
M.\ Liebeck and D.\ Testerman \cite{liebecktesterman}.

\begin{prop}
\label{liebecktester}
Suppose that $H$ is a 
(non-connected) $G$-irreducible  subgroup of $G$ where 
$H^0$ is not $G$-irreducible. Then 
$C_G(H^0)$ contains a non-central torus of $G$.
\end{prop}

\begin{proof}
As $H$ is $G$-ir, it is $G$-cr, so $H^0$ is $G$-cr, by
Theorem \ref{Serrequestion}.  
By hypothesis, $H^0$ is not $G$-ir, 
so there exists a proper parabolic subgroup of $G$ containing $H^0$. 
As $H^0$ is $G$-cr, it lies in a Levi subgroup of this 
parabolic subgroup. This Levi subgroup is  
the centralizer in $G$ of some non-central  torus $S$ of $G$ 
(since the Levi subgroup is proper in $G$).
Thus $C_G(H^0)$ contains $S$, as desired.
\end{proof}

\subsection{$G$-complete reducibility and semisimple modules}
\label{sub:ss}

There is a fundamental connection between Serre's notion of 
$G$-complete reducibility and the semisimplicity of $G$-modules.
Let $h$ denote the Coxeter number of $G$.
For a finite-dimensional $G$-module $V$ define
$n(V)=\max \{\sum_{\alpha>0}\inprod{\lambda,\alpha\check}\}$, 
where the maximum is taken over all $T$-weights $\lambda$ of $V$.
Observe that if $V$ is \emph{non-degenerate}, 
i.e., when the connected kernel of the representation of $G$ on $V$ 
is a torus, then $n(V)\ge h-1$, cf.\  \cite[p.\  20]{serre1}.
The following is the main objective in \cite{serre1}, see also 
\cite[Thm.\ 5.4]{serre2}.

\begin{thm}  
\label{serremain}
Let $V$ be a $G$-module with $p>n(V)$.  Let
$H$ be a closed subgroup of $G$.  
\begin{itemize}
\item[(i)]
If $H$ is $G$-completely reducible, 
then $V$ is a semisimple $H$-module.  
In particular, $V$ is a semisimple $G$-module.
\item[(ii)] If $V$ is non-degenerate and 
semisimple as an $H$-module, then 
$H$ is $G$-completely reducible.
\end{itemize}
\end{thm}

This is a deep theorem and its proof, which is quite involved, 
requires Serre's notion of saturation, 
cf.\ \cite[p.\ 22]{serre1},  \cite[\S 5]{serre2}.
The proof of part (ii) uses the full force of 
Theorem \ref{gcr-thm} below, which itself is
a very difficult case-by-case analysis.

The following is the special case of Theorem \ref{serremain}
for the adjoint module, \cite[Cor.\ 5.5]{serre2}.
Note that $n(\mathfrak g) = 2h-2$ and $\mathfrak g$ 
is a non-degenerate $G$-module.

\begin{cor}
\label{lie-gcr}
Let $H$ be a closed subgroup of $G$.  
Suppose that $p> 2h-2$.
Then $\Lie G$ is a semisimple $H$-module if and only
if $H$ is $G$-completely reducible. 
\end{cor}

\begin{rems}
\label{rems:lie-gcr}
(i). 
For $G = \GL(V)$ the forward implication 
in Corollary \ref{lie-gcr} does not require any 
restrictions on $p$ (see \cite[Thm.\ 3.3]{serre0}; 
this is also a special case of Theorem \ref{thm:adjcr} below,
cf.\ Example \ref{exmp:allsep}). 

(ii). 
It follows from Corollary \ref{red-pair-gl} 
and Example \ref{ex:adjoint} that for 
the forward implication 
in Corollary \ref{lie-gcr}, it suffices to require that
$p > 2$ and $p$ does not divide 
$n+1$, $2n-1$, $n+1$, $n-1$, in case $G$ is an adjoint simple group of  type 
$A_n$, $B_n$, $C_n$, $D_n$, respectively, and 
that $p$ is good for $G$ in case $G$ is an adjoint simple group 
of exceptional type.
For example, if $G$ is adjoint of type $C_n$, then it suffices to require that 
$p > 2$ and $p$ does not divide $n+1$,
which improves on the bound $p>2h-2 = 4n - 2$ from Corollary \ref{lie-gcr}; 
if $G$ is of type $E_8$, then we obtain the bound $p>5$,
which improves on the bound $p>2h-2=58$ from Corollary \ref{lie-gcr}.

(iii). 
It is shown in \cite[Cor.\ 3]{liebeckseitz} 
that if $G$ is simple 
of exceptional type and $H$ is simple of rank at least $2$,
then it suffices to require that $p > 7$ for the reverse implication 
of Corollary \ref{lie-gcr}  to hold.
A computation of J.-P.\ Serre
shows (except possibly when $G = G_2$ and $p=5$ or $G = B_3$ and $p=3$)
that if $p \le 2h - 2$, then there always exists 
a subgroup $H$ of $G$ of type $A_1$ which is $G$-cr, while $\Lie G$ 
is not $H$-semisimple, cf.\ \cite[Rem.\ 5.6]{serre2}.
In particular, this says that the bound in Corollary \ref{lie-gcr} is sharp
for the reverse implication. 
\end{rems}

To go with the counterexamples to the reverse implication of
Corollary \ref{lie-gcr} when  $p \le 2h-2$ 
referred to in Remark \ref{rems:lie-gcr}(iii), 
here is another example of the failure of 
the conclusion of Theorem \ref{serremain}(i) if $p \le n(V)$.

\begin{exmp}
\label{ex:serrecountereg}
Let $p >2$.
The adjoint representation of 
$H := \PGL_p(k)$ on $\hh$ has a 
composition factor of dimension $m:=p^2-2$.
So we can regard $H$ as an $\SL_m(k)$-irreducible 
subgroup of $\SL_m(k)$.  The image of $H$ 
in $G := \PGL_m(k)$ is $G$-irreducible, 
by Lemma \ref{lem:epimorphisms}(ii)(a), 
and it is clear that this image is isomorphic to $H$. So we 
can regard $H$ as a $G$-irreducible subgroup of $G$.  Now $\gg$ is 
simple as a $G$-module, as $p$ is coprime to $m$. 
However,  $\gg$ is not semisimple 
as an $H$-module: for the $H$-submodule $\hh$ of $\gg$ is not 
semisimple. Note that $n(\gg) = 2m-2$.
\end{exmp}

The following example, due to M.W.\ Liebeck,
shows that Theorem \ref{serremain}(ii) fails without the restriction on $p$.

\begin{exmp}
\label{ex:liebeckcountereg}
Suppose that $p=2$, $m\geq 4$ is even and $H := \SP_m(k)$ 
is embedded diagonally in the maximal rank subgroup 
$M := \SP_m(k) \times \SP_m(k)$ of
$G := \SP_{2m}(k)$.  Let $V$ and $V'$ be the natural 
modules for the $\SP_m(k)$-factors of $M$, so 
that the orthogonal direct sum $W :=V\oplus V'$ is the natural module 
for $G$.  Then, as $V$ and $V'$ are irreducible $\SP_m(k)$-modules,  
$W$ is a semisimple $H$-module.  Choose a form-preserving $\SP_m(k)$-module 
isomorphism $f:V\ra V'$.  It is easily checked that the only $m$-dimensional 
$H$-stable subspaces of $W$ are $V'$ and the subspaces $V_a$ 
defined by $V_a:=\{v+af(v)\,|\,v\in V\}$, where $a\in k$.  
Since $(v_1+af(v_1),v_2+af(v_2))= (1+a^2)(v_1,v_2)$ for 
every $v_1,v_2\in V$, we see that $V_1$ is the unique 
$H$-stable, totally singular $m$-dimensional subspace of 
$W$.  This implies that $H$ lies in a proper 
parabolic subgroup $P$ of $G$, but $H$ does not lie 
in any Levi subgroup of $P$.  Thus $H$ is not $G$-cr. 

On the other hand, it follows from Lemma \ref{lem:epimorphisms}(i) 
(applied to the diagonal embedding $H \to M$)
that $H$ is $M$-cr,  
and further from Corollary \ref{cor:linearlyreductive}
that if $p\ne 2$, then $H$ is $G$-cr, as the maximal rank subgroup $M$
is the centralizer of an involution in $G$.
This improves on the bound $p > n(W) = 2m = h$ from 
Theorem \ref{serremain}(ii) in this particular case.
\end{exmp}

Our next result gives another sufficient condition for the 
forward direction of Corollary \ref{lie-gcr} to hold;
recall Definition \ref{def:separable} of a separable subgroup of $G$.

\begin{thm}
\label{thm:adjcr}
Let $H$ be a separable subgroup of $G$. 
If $\gg$ is semisimple as an $H$-module, 
then $H$ is $G$-completely reducible.
\end{thm}

\begin{proof}
By Remark \ref{rem:topfg},
we can assume that $H$ is topologically 
finitely generated, say by $h_1,\ldots, h_n$.
Suppose that $H$ is not $G$-cr.  Then, by Corollary \ref{cor:closedorbit},
the $G$-orbit of $(h_1,\ldots, h_n)$ is not closed in $G^n$.  
By the Hilbert--Mumford Theorem \ref{thm:HMT} this implies 
that there exists a cocharacter $\lambda$ of $G$ such that 
$\underset{t\to 0}{\lim}\,\lambda(t)\cdot (h_1,\ldots, h_n)$ exists ---
call this limit $(h_1',\ldots, h_n')$ ---  and such that 
$G\cdot (h_1',\ldots, h_n')$ is closed.  In particular, we have
\begin{equation}
\label{ineq}
\dim (G\cdot (h_1',\ldots, h_n'))  < \dim (G\cdot (h_1,\ldots, h_n)).
\end{equation}
Let $H'$ be the algebraic group generated by $h_1',\ldots, h_n'$, 
so that $H'=c_\lambda(H)$, where $c_\lambda$ is the map defined in 
Lemma \ref{lem:cochars}(iii). Note that 
$H \subseteq P_\lambda$, since 
$\underset{t\to 0}{\lim}\,\lambda(t)\cdot h_i$ exists
for each $i = 1, \ldots, n$.
Inequality \eqref{ineq} implies that 
$\dim C_G(H')  > \dim C_G(H)$.  Since $H$ is separable in $G$, we deduce 
that $\dim \cc_\gg(H')  > \dim \cc_\gg(H)$.
 
Let $M=\Ad_G(H)$ and $M'=\Ad_G(H')$.  
It is clear that $M' = c_{\Ad_G \circ \lambda}(M)$.  
Since $\gg$ is $H$-semisimple, 
$M$ is $\GL(\gg)$-cr. It follows from Lemma \ref{lem:srcrit} 
and Theorem \ref{mainthm}
that $M'$ is $\GL(\gg)$-conjugate to $M$.  Now $\cc_\gg(H)$ 
(respectively $\cc_\gg(H')$) 
is precisely the set of fixed points of $M$ (respectively $M'$) 
in $\gg$, so $\cc_\gg(H')$ 
is $\GL(\gg)$-conjugate to $\cc_\gg(H)$.  But this implies that 
$\dim \cc_\gg(H')  = \dim \cc_\gg(H)$, a contradiction.  
We conclude that $H$ is $G$-cr, as required. 
\end{proof}

The following is a simplified statement of 
a consequence of a number of deep theorems
due to J.C.~Jantzen \cite{Jantzen0} and G.~McNinch  \cite{mcninch}
in case $G$ is classical 
and M.~Liebeck and G.~Seitz \cite{liebeckseitz0}
for $G$ of exceptional type;
see \cite[Thm.\ 4.4]{serre2} and \cite[\S 3]{serre1.5}.

\begin{thm}
\label{gcr-thm}
Let $H$ be a closed connected subgroup of $G$ and suppose that $p > h$.
Then $H$ is $G$-completely reducible if and only if
$H$ is reductive. 
\end{thm}

Theorem \ref{gcr-thm} says that provided $p > h$, 
for a connected subgroup of $G$
the notions of reductivity and $G$-complete reducibility 
are equivalent, as in characteristic zero, cf.\ Subsection \ref{sub:charzero}. 
As the results in \cite{liebeckseitz0} and  \cite{mcninch}
depend on case-by-case studies,
it would be desirable to have a uniform proof of 
Theorems \ref{serremain} and \ref{gcr-thm}, 
even with some additional restrictions on $p$.
We believe that Theorem \ref{red-pair}, Corollary \ref{red-pair-gl}, 
and Theorem \ref{thm:adjcr}
provide a first step in this direction.  

If $(\GL(V),H)$ is a reductive pair, then $V$ is non-degenerate by definition. 
The dependence on the characteristic in 
Corollary \ref{red-pair-gl} 
is buried in the hypothesis that $(\GL(V),H)$ is a reductive pair, see
Example \ref{exmp:redpair},
Remark \ref{rem:redpair}, and Example \ref{ex:adjoint}.
The advantages of Theorem \ref{red-pair},  
Corollary \ref{red-pair-gl}, and Theorem \ref{thm:adjcr} 
are that they do not require the notion of saturation and 
they are free of case-by-case considerations.

\section{$G$-Complete Reducibility and Buildings}
\label{s:building}

In this section we consider 
the connection between the notion of 
$G$-complete reducibility and the building of $G$ due to 
J.-P.\  Serre \cite[Thm.\ 2]{serre1.5}. 
For an arbitrary spherical building $X$, a subset $Y$ of
$X$ is said to be \emph{convex} if whenever two points of  
$Y$ are not opposite in $X$, then  $Y$  
contains the unique geodesic joining these points.
A convex subset $Y$ is 
\emph{$X$-completely reducible} ($X$-cr) if
for every $y\in Y$ there
exists a point $y'\in Y$ opposite to $y$ in $X$, 
\cite[Def. 2.2.1]{serre2}.
The vertices of $X$ can be labelled in an 
essentially unique way via an equivalence
relation on vertices, cf. 
\cite[p30]{brown}, \cite[2.1.2]{serre2}; 
the \emph{type} of a vertex is its label.
An automorphism $f$ of $X$ is said to be \emph{type-preserving}
if $x$ and $f(x)$ have the same type for all vertices $x$ of $X$.

Now let $X = X(G)$ be the spherical Tits building of $G$, 
cf.\ \cite{brown}, \cite{tits1}.
Recall that the simplices in $X$ correspond to
the parabolic subgroups of $G$ and the vertices
of $X$ correspond to the maximal proper parabolic subgroups, see
\cite[\S 3.1]{serre2}.
For a subgroup $H$ of $G$ 
let $X^H$ be the fixed point subcomplex of the action of $H$, i.e.,
the subcomplex of all $H$-stable (thus $H$-fixed) simplices in $X$.
This subcomplex is always convex; if it
is also $X$-completely reducible, then we say 
\emph{$H$ acts completely reducibly on $X$}, \cite[\S  2.3]{serre2}.
For any subgroup $H$ of $G$ the action of $H$ on $X$ is type-preserving.

Take a Levi subgroup $L$ of $G$ and let $s(L) := X^L$ denote the 
subcomplex of $X$ consisting of the parabolic subgroups of $G$
containing $L$.
For every parabolic subgroup $P$ in $s(L)$ there is a unique 
Levi subgroup $M$ of $P$ with $L\subseteq M$.
Moreover, the parabolic subgroup $P^-$  
such that $P\cap P^- = M$ is also contained in $s(L)$, so that
$P$ has an opposite in $s(L)$;
thus $s(L)$ is $X$-cr.
This argument also shows that each $P$ in $s(L)$ has a unique opposite
in $s(L)$, and this implies that the
geometric realization of $s(L)$ has the homotopy type of a 
single sphere (cf.\ Theorem \ref{gcr-building}(v) below). 
Serre calls the subcomplexes $s(L)$ of $X$ \emph{Levi spheres}, 
\cite[\S 2]{serre1.5} or \cite[2.1.6, 3.1.7]{serre2}.
The following is part of \cite[Thm.\ 2]{serre1.5} 
in our context.

\begin{thm}
\label{gcr-building}
Let $H$ be a closed subgroup of $G$. 
Then the following 
are equivalent:
\begin{itemize}
\item[(i)] $H$ is $G$-completely reducible;
\item[(ii)] $X^H$ is $X$-completely reducible;
\item[(iii)] $X^H$ contains a Levi sphere of the same dimension as $X^H$; 
\item[(iv)] $X^H$ is not contractible 
(i.e., does not have the homotopy type of a point);
\item[(v)] 
$X^H$ has the homotopy type of a bouquet of spheres.
\end{itemize}
\end{thm}

J.-P.\  Serre observed
that Theorem \ref{mainthm} can also be interpreted in 
terms of Theorem \ref{gcr-building}.
Note that any Levi sphere $s(L)$ contained in  $X^H$
is of the form $s(C_G(T))$, where $T$ is a subtorus of $C_G(H)$.
Thus  $s(L)$ 
is of maximal dimension in  $X^H$ if and only if $L = C_G(S)$, where 
$S$ is a maximal torus  of $C_G(H)$.
Moreover, $\dim X^H  =  \dim s(L)$  if and only if
$H$  is not contained in a proper parabolic subgroup of the Levi
subgroup  $L$. 
Thus part (iii) of  Theorem \ref{gcr-building}
is equivalent to $H$ being strongly reductive in $G$.
We are grateful to J.-P.\  Serre for this 
building-theoretic interpretation. 

Thanks to Theorem \ref{gcr-building}, the  
results of the previous section have counterparts in terms
of buildings; e.g., Corollary \ref{cor1converse}
then says that for a closed subgroup $H$ of $G$
the fixed point subcomplex $X^H$ 
is contractible if and only if $X^{N_G(H)}$ is contractible, etc.

If $H$ is a $G$-cr subgroup of $G$, then it follows from 
Theorem \ref{gcr-building} that $X^H$ is itself a spherical chamber 
complex, much like a building.  
However, in general $X^H$  is not a building.
For instance, let $T$ be a maximal torus in $G$.
Then $T$ is $G$-cr. However, $X^T$ is a Coxeter complex and so is not 
a building, as  $X^T$ is not thick.
(Recall that a chamber complex is said to be \emph{thick}
if every simplex of codimension one is contained 
in at least three chambers, e.g., see \cite{brown}.)
See also Example \ref{ex:f4} below.  

The following lemma extends a result that is well known when $H$ is a torus.

\begin{lem}
\label{parabolics-cent}
Let $H$ be a $G$-completely reducible subgroup of $G$ and
let $P$ be a parabolic subgroup of $G$ containing $H$.
Then $P\cap C_G(H)^0$ is a parabolic subgroup of $C_G(H)^0$.
Moreover, all parabolic subgroups of $C_G(H)^0$ arise in this way.
\end{lem}

\begin{proof}
First notice that since $H$ is $G$-cr, $C_G(H)^0$
is reductive by Proposition \ref{redcent}.
Since $H$ is $G$-cr, there exists $\lambda \in Y(G)$ such that 
$H \subseteq L_\lambda(G) \subseteq P_\lambda(G) = P$ by 
Lemma \ref{lem:cochars}.
Then $\lambda(k^*) \subseteq C_G(H)^0$, so that $\lambda \in Y(C_G(H)^0)$ and
$P\cap C_G(H)^0 = P_\lambda(C_G(H)^0)$ 
is a parabolic subgroup of $C_G(H)^0$, by Lemma \ref{lem:cochars},
as claimed.

Now for any parabolic subgroup $Q$ of $C_G(H)^0$,
there exists $\mu \in Y(C_G(H)^0)$ such that $Q = P_\mu(C_G(H)^0)$. 
Then the parabolic subgroup $P = P_\mu(G)$ of $G$
contains $H$, since $\mu(k^*)$ centralizes $H$,
and $Q = P \cap C_G(H)^0$ by Corollary \ref{cor:parsofsubgp}.
\end{proof}

For an arbitrary reductive group $K$, let $X(K)$ denote
the spherical Tits building of $K$.
In our next result, Lemma \ref{parabolics-cent} allows us to relate 
the fixed point complex $X(HC_G(H)^0)^H$
and the building of the connected centralizer of $H$ in
case $H$ is a connected $G$-completely reducible subgroup of $G$.  

\begin{prop}
\label{fixbuilding-regular}
Let $H$ be a connected $G$-completely reducible subgroup of $G$.
Then the chamber complexes
$X(HC_G(H)^0)^H$ and $X(C_G(H)^0)$ are isomorphic.
In particular, we see that $X(HC_G(H)^0)^H$ is itself a building.
\end{prop}

\begin{proof}
Set $M = HC_G(H)^0$ and observe that $M$ is reductive by 
Proposition \ref{redcent}.
Since $H$ is normal in $M$, $H$ is $M$-cr by
Theorem \ref{Serrequestion}.
Also, since $C_G(H)^0 \subseteq M$, we have $C_M(H)^0 = C_G(H)^0$.

By Lemma \ref{parabolics-cent} above, parabolic
subgroups $P$ of $M$ containing $H$ 
correspond bijectively to parabolic 
subgroups $P\cap C_M(H)^0$ of $C_M(H)^0 = C_G(H)^0$:
to see that the correspondence is one-to-one, notice that if 
$H\subseteq P \subseteq M = HC_G(H)^0$, then we have
$P = H(P\cap C_G(H)^0)$. 
Thus, for two parabolic subgroups $P$ and $Q$ of $M$ containing $H$, we
see that $P = Q$ if and only if $P\cap C_G(H)^0 = Q\cap C_G(H)^0$. 

This bijection on the sets of parabolic subgroups affords the desired
isomorphism on the underlying chamber complexes.
\end{proof}

To illustrate Proposition \ref{fixbuilding-regular} we first point 
to a trivial case:

\begin{exmp}
\label{ex:normal}
Let $G$ be a reductive group and let $H$ be a proper closed 
connected normal 
subgroup of $G$. Since $H$ is normal in $G$, $H$ is $G$-cr, 
by Theorem \ref{Serrequestion}. 
As a special case of Proposition \ref{fixbuilding-regular} we obtain
that $X^H$ is isomorphic to $X(C_G(H)^0)$.
\end{exmp}

The following  example was communicated to 
us by B.\ M\"uhlherr.

\begin{exmp}
\label{ex:f4}
Let $G$ be of type $F_4$ and let $H$ be the simple $A_1$-factor
of the maximal rank subgroup $M$ of $G$ of type $A_1C_3$.
Then, by Proposition \ref{regular-is-cr}, $H$ is $G$-cr, so that 
$X^H$ does have the 
homotopy type of a bouquet of spheres, by Theorem \ref{gcr-building}.
However, $X^H$  is not a building, as the apartments in $X^H$ 
are not Coxeter complexes.
Note that we have $HC_G(H) = M$. 
It follows from Proposition \ref{fixbuilding-regular} that
$X(M)^H$ is a building of type $C_3$.
\end{exmp}

Group actions on buildings 
are also considered in \cite{muhlherr} and \cite{muhlherrschmid}
(Section 5 in \cite{muhlherrschmid} is closely related to 
Theorem \ref{gcr-building}).
Given a spherical chamber complex, one can associate to 
it a building via a so-called \emph{thickening} procedure,
see \cite[\S  1.7]{muhlherr} and \cite[\S 5]{muhlherrschmid} for details.
It follows from \cite[1.7.26, 1.8.22, 3.4.8]{muhlherr} that
for a $G$-completely reducible 
subgroup $H$ of $G$, the fixed point subcomplex $X^H$ 
is a building if and only if it is thick.
Moreover, for any $G$-completely reducible
subgroup $H$ of $G$, 
the thickening of $X^H$ is isomorphic to 
$X(C_G(H)^0)$.
Note that if $X^H$ is already a building, then thickening
has no effect.
These results hold in greater generality than is stated
here.
For example, they are true for any group $H$ 
which acts on $G$ such that the induced action on $X$ 
is completely reducible and type-preserving,
see \cite{muhlherr}. 

\smallskip

Now suppose that $Y$ is a \emph{strictly convex} subcomplex of $X$, 
i.e., $Y$ is convex but it
does not contain any two opposite points of $X$.
Suppose the subgroup $H$ of $G$ stabilizes $Y$.
The so-called ``Center Conjecture'' due to J.\  Tits 
claims the existence of a 
fixed point of $H$ in $Y$, cf.\  \cite[\S 4]{serre1.5}. 
It turns out that Theorem \ref{Serrequestion} is a 
consequence of Tits' Center Conjecture, \cite[Prop.\ 2.11]{serre2}.
We refer to \cite[\S 2.4]{serre2} for more details and known instances of 
this conjecture; see also \cite[\S 3.6]{muhlherr} and
\cite{muhlherr2} for related results.

\section{Rationality Questions}
\label{s:rationality}

The notions of $G$-complete reducibility, etc., can be extended 
to reductive groups defined over arbitrary fields; see \cite{serre2}.  
In this section $k$ denotes an arbitrary field, not necessarily 
algebraically closed, $\overline{k}$ denotes the algebraic closure 
of $k$, and $G$ denotes a reductive 
group defined over $k$ (see \cite{borel}, \cite{boreltits}, and \cite{spr2}
for more details).  
Given a field extension $k'/k$, 
we denote by $G(k')$ the group 
of $k'$-rational points of $G$; 
if $k'$ is algebraically closed, 
then we often identify $G$ with $G(k')$.  
By a $k'$-subgroup $H$ of $G$, we mean an algebraic
subgroup of $G$ over $k'$.
We call a parabolic subgroup of $G$ that is defined 
over $k'$ a $k'$-parabolic subgroup of $G$.  If $P$ is a 
$k'$-parabolic subgroup of $G$ and $L$ is a Levi 
subgroup of $P$ that is defined over $k'$, then we 
call $L$ a $k'$-Levi subgroup of $P$. 

\begin{defn}
 Let $k'/k$ be a field extension.  We say that a $k'$-subgroup $H$ of $G$ is 
\emph{$G$-completely reducible over $k'$} if whenever 
$H$ is contained in a $k'$-parabolic subgroup $P$ of 
$G$, there is a $k'$-Levi subgroup $L$ of $P$ 
such that $H\subseteq L$.
\end{defn}

\begin{rem}
 In particular, if $k'$ is algebraically closed, 
then $H$ is $G$-completely reducible over $k'$ 
if and only if $H(k')$ is $G(k')$-completely 
reducible as defined in Section \ref{s:intro}.
\end{rem}

\subsection{Algebraically closed fields}
\label{sub:algclos}

First we consider extensions $k'/k$ of algebraically closed fields.
In view of Theorem \ref{mainthm}, we record some results 
from \cite{martin1}, replacing ``strongly reductive in $G$'' 
with ``$G$-cr'':
parts (i) and (ii) of the following theorem are 
\cite[Prop.\ 10.2]{martin1} and
\cite[Thm.\ 10.3]{martin1}, respectively.

\begin{thm}
\label{thm:algclosrat}
Let $k'/k$ be an extension of algebraically closed fields.
\begin{itemize}
\item[(i)] 
Let $H$ be a $k$-subgroup of $G$.
Then $H$ is $G$-completely reducible over $k'$ if and only if $H$ 
is $G$-completely reducible over $k$.
\item[(ii)] 
Let $K$ be a $k'$-subgroup of $G$ such that 
$K$ is $G$-completely reducible over $k'$.  
Then there exists a $k$-subgroup $H$ of $G$ 
such that $H$ is $G$-completely reducible over $k$ 
and $H$ is $G(k')$-conjugate to $K$.
\end{itemize}
\end{thm}

\begin{rem}
\label{rem:algclosratirind}
 Part (i) of Theorem \ref{thm:algclosrat} also holds replacing 
``completely reducible'' with ``irreducible'' or ``indecomposable'': 
for a closed subgroup is $G$-ind if and only if it is not centralized 
by any non-central torus of $G$ and is $G$-ir if and only if it is both 
$G$-cr and $G$-ind.
\end{rem}

If $k'/k$ is an extension of algebraically closed fields, 
then two $k$-subgroups $H_1$ and $H_2$ of $G$ are 
$G(k')$-conjugate if and only if they are $G(k)$-conjugate.  
This, together with Theorem \ref{thm:algclosrat}, 
establishes the following result.
\begin{cor}
\label{cor:algexteqvc}
Let $k'/k$ be an extension of algebraically closed fields.
Then the map $H(k)\mapsto H(k')$ gives rise to a bijection
between the set of $G(k)$-conjugacy classes of 
$k$-subgroups of $G$ that are $G$-completely reducible over $k$, 
and the set of $G(k')$-conjugacy classes of 
$k'$-subgroups of $G$ that are $G$-completely reducible over $k'$.
\end{cor}

\begin{cor}
\label{cor:countcr}
Suppose that $k$ is algebraically closed.
Then there are only countably many conjugacy 
classes of $k$-subgroups of $G$ that are $G$-completely reducible over $k$.
\end{cor}

\begin{proof}
Let $k_0$ be the algebraic closure of the prime field of $k$ and let 
$k_1$ 
be the algebraic closure of $k_0(t)$, where $t$ is transcendental 
over $k_0$.  The group $G$ admits a $k_0$-structure by 
\cite[Prop.\ 3.2]{martin1}, so we can assume that $G$ 
is defined over $k_0$.  By 
Corollary \ref{cor:algexteqvc}, we can assume that 
$k=k_1$; in particular, $k$ is countable and $k/k_0$ 
is transcendental.  Now any $k$-subgroup of $G$ that 
is $G$-cr over $k$ is reductive, and hence by \cite[Lem.\ 9.2]{martin1}
is topologically finitely 
generated.  But $G(k)$, being countable, has only 
countably many topologically finitely generated 
subgroups, so the result follows.
\end{proof}

\begin{rem}
We cannot replace ``countably many'' by ``finitely many'' in the previous 
corollary: 
see Remark \ref{rem:conjclass}(ii).
\end{rem}

\subsection{Perfect fields}
\label{sub:perfect}

Now we consider field extensions $k'/k$ where both $k'$ and $k$ are perfect.  
If $k$ is perfect, $k'/k$ is a Galois extension, and $X$ is a variety defined 
over $k$, then a $k'$-subvariety of $X$ is defined over $k$ 
if and only if it is $\Gal(k'/k)$-stable,
cf.\ \cite[AG \S 14]{borel}.  The forward implication 
of the following result
uses the argument of \cite[Prop.\ 2.2]{LMS}, but 
without the extra complication 
that appears there of graph automorphisms;
\cite[Prop.\ 2.2]{LMS} relies on a version of the Hilbert--Mumford--Kempf 
Theorem (see \cite[Thm.\ 4.2]{kempf}) which requires 
$k$ to be perfect.  The reverse implication of Theorem \ref{thm:perfectcr}
is an observation of J.-P.\ Serre.

\begin{thm}
\label{thm:perfectcr}
Let $k'/k$ be an extension of 
perfect fields and let $H$ be a $k$-subgroup of $G$.
Then $H$ is $G$-completely reducible over $k'$ if and only if $H$ is 
$G$-completely reducible over $k$.
\end{thm}

\begin{proof}
It suffices to prove the result in the special 
case when $k'$ is algebraically closed.  
By Theorem \ref{thm:algclosrat}(i), we can 
assume that $k'= \overline{k}$.  Suppose that 
$H$ is not $G$-cr over $\overline{k}$.
By Lemma \ref{lem:fgapprox}, there exists a finitely 
generated subgroup $\Gamma$ of $H(\overline{k})$ such 
that for every $\overline{k}$-parabolic 
subgroup $P$ of $G$ and every $\overline{k}$-Levi 
subgroup $L$ of a $\overline{k}$-parabolic subgroup of $G$, we have 
$\Gamma\subseteq P(\overline{k})$ if and only if $H\subseteq P$ and 
$\Gamma\subseteq L(\overline{k})$ 
if and only if $H\subseteq L$.
 
As $k$ is separable and $\overline{k}/k$ is algebraic, 
we can choose a finite Galois 
extension $k_1/k$ such that $\Gamma\subseteq H(k_1)$.  
Let $\Gamma_1$ be the group generated by the 
$\Gal(\overline{k}/k)$-conjugates of $\Gamma$ and 
let $M$ be the closure of $\Gamma_1$ in $H$.  Then $\Gamma_1$ is a 
finitely generated $\Gal(\overline{k}/k)$-stable 
subgroup of $G(\overline{k})$, and we can 
choose a finite set of generators $h_1,\ldots, h_n$ for $\Gamma_1$ such 
that the $h_i$ are permuted by $\Gal(\overline{k}/k)$.  
By Theorem \ref{mainthm} and the argument of 
\cite{LMS} applied to the tuple $(h_1,\ldots, h_n)$ (see the proof 
of \cite[Prop.\ 2.2]{LMS} and the paragraph that follows it), there 
exists a $\Gal(\overline{k}/k)$-stable $\overline{k}$-parabolic 
subgroup $P$ of $G$ such 
that $M\subseteq P$ but $M$ does not lie in any $\overline{k}$-Levi 
subgroup of $P$.  
Then $P$ is defined over $k$, $H\subseteq P$ and $H$ does not 
lie in any $\overline{k}$-Levi subgroup of $P$; in particular, $H$ does not 
lie in any $k$-Levi subgroup of $P$.  It follows that $H$ is not 
$G$-cr over $k$. 
  
Conversely, suppose that $H$ is $G$-cr over $\overline{k}$.  Let $Q$ be a 
$k$-parabolic 
subgroup of $G$ such that $H\subseteq Q$.  
The centralizer $C_Q(H)$ is defined over $k$, 
because $H$ and $Q$ are, so  
by \cite[Thm.\  18.2]{borel},
$C_Q(H)$ contains a maximal torus $S$ defined over 
$k$.  As $H$ is $G$-cr over $\overline{k}$, $H$ is contained 
in a $\overline{k}$-Levi subgroup $M$ of $Q$.  Conjugating $M$ 
by some element of 
$C_Q(H)$ if necessary, we can assume that the torus $Z(M)^0$ 
is contained in $S$.  We have 
$C_G(S)\subseteq C_G(Z(M)^0)=M$.  Now $C_G(S)$ 
is defined over $k$, 
because $S$ is, so $C_G(S)$ contains a maximal torus $T$ defined 
over $k$, again by \cite[Thm.\  18.2]{borel}.
There is exactly one Levi subgroup of $Q$ containing any given 
maximal torus, cf.\ \cite[Cor.\ 8.4.4]{spr2},
so it follows that $M$ is defined over $k$.
Thus $H$ is $G$-cr over $k$.
\end{proof}

\begin{rem}
One can extend the definition of $G$-ir and $G$-ind to the 
non-algebraically closed setting in the obvious way, but the 
analogue of Theorem \ref{thm:perfectcr} does not carry over 
(cf.\ Remark \ref{rem:algclosratirind}); 
for example, take $k$ perfect but not 
algebraically closed and consider a $k$-subgroup $H$ of $G = \GL_n(k)$ that 
is irreducible but not absolutely irreducible.
\end{rem}

The following special case has applications to 
finite groups of Lie type (cf.\ \cite[Prop.\ 2.2]{LMS}). 

\begin{exmp}
\label{exmp:Lietype}
Suppose that $G$ is defined over the finite 
field ${\mathbb F}_p$, let $k'$ be the algebraic 
closure of ${\mathbb F}_p$, and let $\sigma : G(k')\ra G(k')$ 
be some power of Frobenius.  Let $G_\sigma$ be the subgroup 
of fixed points of $\sigma$.
Then any subgroup 
$F$ of $G_\sigma$ either is $G$-completely reducible over $k'$, or 
is contained in a proper 
$\sigma$-stable parabolic subgroup of $G$.  
To see this, observe that we have $G_\sigma=G(k)$ 
for some finite extension $k$ of 
${\mathbb F}_p$.  We can regard $F$ as a $k$-subgroup of $G$.  
If $F$ is not $G$-cr over $k'$, then $F$ is not $G$-cr over $k$, 
by Theorem \ref{thm:perfectcr}, 
so there is a $k$-parabolic subgroup $P$ of $G$ such 
that $F$ is not contained in any $k$-Levi subgroup of $P$.  
In particular, $P$ is a proper $\sigma$-stable parabolic 
subgroup of $G$, as required.
If $P$ is chosen as in the proof of Theorem \ref{thm:perfectcr}, 
then in fact $F$ is not contained in any $k'$-Levi subgroup of $P$ at all.
\end{exmp} 

We are grateful to G.~McNinch for the following example, 
which shows that the reverse implication in Theorem \ref{thm:perfectcr} 
fails if the hypothesis of perfection is removed.

\begin{exmp}
\label{exmp:nonperfect}
Let $k_1/k$ be a purely inseparable field extension 
of degree $p$.  Set $k'= \overline{k}$.  We can regard 
$k_1^*$ as an algebraic group $H$ over $k$; the action 
of $k_1^*$ on $k_1$ by left multiplication gives rise 
to a $k$-embedding of $H$ in $G = \GL(k_1)$, where we regard 
$k_1$ as a $k$-vector space.  As $H$ acts transitively on 
$k_1^*$, $H$ cannot stabilize any proper non-trivial $k$-subspace 
of $k_1$, so $H$ is $G$-cr over $k$.  However, $H$ is not reductive: 
for the homomorphism $\phi$ sending $x$ to $x^p$ maps $H$ onto the group 
of scalar multiples of the identity matrix, so $\ker \phi$ is a 
$(p-1)$-dimensional normal unipotent subgroup of $H$.  This implies 
that $H$ is not $G$-cr over $k'$.
\end{exmp}

\section{The Non-Connected Case}
\label{s:noncon}

In this section we extend Theorem \ref{mainthm} 
and many of our other results to groups $G$ 
that are no longer required to be connected.  The formalism we 
use here for dealing with such groups is taken from \cite{martin1}, 
which in turn is based on the approaches of Vinberg \cite{vinberg} 
and Richardson \cite{rich1}.  
Thus for the remainder of 
the paper we suppose that $G$ is a linear algebraic group with $G^0$ 
reductive; we call this ``the non-connected case'' and
such a group is referred to as
a ``non-connected reductive group''; note, however, that we
do not exclude the case $G=G^0$.

The idea is to use the appropriate generalization of the 
notion of a parabolic subgroup to the non-connected case, 
using the formalism of Lemma \ref{lem:cochars} (cf.\ \cite[\S 5]{martin1}).
For $\lambda\in Y(G)$, we call a subgroup of the form 
$P_\lambda := \{g\in G \mid 
\underset{t\to 0}{\lim}\, \lambda(t) g \lambda(t)\inverse\ \textrm{exists}\}$ 
a {\em Richardson parabolic} 
(or {\em R-parabolic}) subgroup of $G$ (cf.\ \cite[\S 2]{rich1}; 
in \cite{martin1}, these were called generalized parabolic subgroups).  
By \cite[Lem.\  6.2.4]{spr2}, any R-parabolic subgroup $P$ of $G$ 
is a parabolic subgroup of $G$ in the sense that $G/P$ is a complete variety, 
but the converse is false,  cf.\ \cite[Rem.\ 5.3]{martin1}.    
We call a subgroup of the form $L_\lambda:=C_G(\lambda(k^*))$ a 
{\em Richardson Levi} (or {\em R-Levi}) subgroup of $P_\lambda$.  
By an R-Levi subgroup of $G$, we mean an R-Levi subgroup of some 
R-parabolic subgroup of $G$.  
We have $P_\lambda=L_\lambda\ltimes R_u(P_\lambda)$.  
If $L$ is an R-Levi subgroup of $G$, then $L=C_G(Z(L)^0)$
and $C_G(L) = Z(L)$.

We define $G$-complete reducibility, $G$-irreducibility and 
$G$-indecomposability as we did in the connected case, 
replacing parabolic subgroups and Levi subgroups with 
R-parabolic subgroups and R-Levi subgroups respectively, 
and we extend the definition of strong reductivity in $G$ similarly.  
If $\lambda\in Y(G)$, then $P_\lambda^0=P_\lambda(G^0)$ is a parabolic 
subgroup of $G^0$ and $L_\lambda^0=L_\lambda(G^0)$ is a Levi subgroup 
of $P_\lambda^0$.  
Conversely, if $P$ is any parabolic subgroup of $G^0$, 
then $P=P_\lambda^0$ for some $\lambda\in Y(G)$; 
moreover, if $L$ is a Levi subgroup of $P$, then 
$\lambda$ can be chosen so that $L = L_\lambda^0$. 
Thus the above definitions of $G$-complete reducibility, etc., 
agree with those from Section \ref{s:intro} if $G$ is connected.

The proof of the non-connected version of Theorem \ref{mainthm} 
is given in Subsection \ref{sub:cr}.  We use the theorem to prove 
several of the results (e.g., the non-connected version of 
Proposition \ref{Rich-linearlyreductive}) in 
Subsections \ref{sub:prelim} and \ref{sub:git-gcr}; 
the proof of the non-connected version of 
Theorem \ref{mainthm} does not depend 
on these results.

\subsection{Preliminaries}
\label{sub:prelim}
In order to generalize our work to the 
non-connected case, we need to prove the analogues 
of a number of results which are standard
when $G$ is connected.
Chief among these is \cite[Prop.\ 4.4]{boreltits}
which is central to the proof of Theorem \ref{mainthm}. 

First notice that all of
Lemma \ref{lem:cochars} extends to the non-connected case, 
replacing parabolic subgroups and Levi subgroups with R-parabolic
subgroups and R-Levi subgroups respectively.  
It is clear that $R_u(P_\lambda)=\ker{c_\lambda}$ is connected, 
so $R_u(P_\lambda)=R_u(P_\lambda^0)$; 
in particular,
$L_\lambda$ meets every component of $P_\lambda$.  
For more details, see \cite[\S 5]{martin1}.

The next result is \cite[Prop.\ 5.4(a)]{martin1}.

\begin{prop}
\label{prop:normpar}
Let $P$ be a parabolic subgroup of $G^0$.  
Then $N_G(P)$ is an R-parabolic subgroup of $G$.
\end{prop}

Note that $N_G(P)^0= P$, 
since a parabolic subgroup of a connected reductive 
group is self-normalizing.

Part (ii) of the next lemma provides the 
extension of \cite[Prop.\ 4.4(c)]{boreltits}.
Notice that part (iii) is also standard for connected $G$.

\begin{lem}
\label{lem:nonconn}
\begin{itemize}
\item[(i)]  Let $\lambda,\mu\in Y(G)$ such that $\lambda(k^*)$ 
and $\mu(k^*)$ commute.  
Then for all sufficiently large $m\in {\mathbb N}$, 
we have $P_{m\lambda+\mu}= P_\mu(L_\lambda)\ltimes R_u(P_\lambda)$ and $L_{m\lambda+\mu}=L_\mu(L_\lambda)$.  
In particular, 
$P_{m\lambda+\mu}\subseteq P_\lambda$.
\item[(ii)] If $P$ is an R-parabolic subgroup of $G$ and 
$L$ is an R-Levi subgroup of $P$, then the R-parabolic subgroups of $G$ 
contained in $P$ are precisely the subgroups of the form 
$P'\ltimes R_u(P)$ with $P'$ an R-parabolic subgroup of $L$.
\item[(iii)] If $P$ and $Q$ are R-parabolic subgroups of $G$ 
with R-Levi subgroups $L$ and $M$ respectively, such that 
$L\cap M$ contains a maximal torus $T$ of $G$, 
then 
\begin{equation}
\label{eq:intersection}
P\cap Q=(L\cap M)(L\cap R_u(Q))(R_u(P)\cap M)(R_u(P)\cap R_u(Q)),
\end{equation} 
and $R_u(P\cap Q)$ is the product of the last three factors.
\end{itemize}
\end{lem}

\begin{proof}
(i). The inclusions $P_{m\lambda+\mu}\subseteq P_\lambda$ 
and $L_{m\lambda+\mu}\subseteq L_\lambda$
for large $m$ follow from the proof of 
\cite[Prop.\ 6.7]{martin1}; the second 
inclusion gives 
$L_{m\lambda+\mu}= L_{m\lambda+\mu}(L_\lambda)= L_\mu(L_\lambda)$.  
Let $T$ be a maximal torus of $P_\lambda$ such that $\lambda,\mu\in Y(T)$.  
For sufficiently large $m$, 
we have $\langle m\lambda+\mu,\alpha\rangle > 0$ for all 
$\alpha\in \Psi(R_u(P_\lambda),T)$, 
whence $R_u(P_{m\lambda+\mu})\supseteq R_u(P_\lambda)$.  
Since $P_\lambda=L_\lambda\ltimes R_u(P_\lambda)$, 
it follows that 
$P_{m\lambda+\mu}= (P_{m\lambda+\mu}\cap L_\lambda)\ltimes R_u(P_\lambda)
= P_{m\lambda+\mu}(L_\lambda)\ltimes R_u(P_\lambda)
= P_\mu(L_\lambda)\ltimes R_u(P_\lambda)$. 

(ii). We can write $P=P_\lambda$ and $L=L_\lambda$ for some $\lambda\in Y(G)$.
If $\mu\in Y(L)$, then $\lambda(k^*)$ and $\mu(k^*)$ commute, 
so $P_\mu(L)\ltimes R_u(P)$ is an R-parabolic subgroup of $G$ by part (i).  
Conversely, let $\mu\in Y(G)$ such that $P_\mu\subseteq P$.  
Let $T$ be a maximal torus of $L$.  
If $P_\mu$ is of the form $P'\ltimes R_u(P)$ for some R-parabolic subgroup 
$P'$ of $L$, then any $P$-conjugate of $P_\mu$ is also of this form; 
thus, since any two maximal tori of 
$P$ are $P$-conjugate, 
we can assume that $\mu(k^*)\subseteq T$.  
Since $P_\mu\subseteq P$, we have $P_\mu^0\subseteq P^0$, 
whence $R_u(P)=R_u(P^0)\subseteq R_u(P_\mu^0)=R_u(P_\mu)$ 
(the middle inclusion is a standard result for connected groups).  
Thus 
$P_\mu = (P_\mu\cap L)\ltimes R_u(P)
= P_\mu(L)\ltimes R_u(P)$.

(iii). Write $P=P_\lambda$, $L=L_\lambda$, $Q=P_\mu$, and $M=L_\mu$, 
where $\lambda,\mu\in Y(G)$.  
Since $T$ is a maximal torus of $L$ and $\lambda(k^*)\subseteq Z(L)^0$, 
we have $\lambda(k^*)\subseteq T$, whence $\lambda(k^*)\subseteq P\cap Q$.  
Since $P\cap Q$ is closed, it follows that $P\cap Q$ is $c_\lambda$-stable, 
so we have $P\cap Q= (L\cap Q)(R_u(P)\cap Q)$.  
By a similar argument, $L\cap Q=(L\cap M)(L\cap R_u(Q))$ 
and $R_u(P)\cap Q= (R_u(P)\cap M)(R_u(P)\cap R_u(Q))$, 
and the product decomposition \eqref{eq:intersection}
of $P\cap Q$ follows.  
It is easily checked that 
\[V:=(L\cap R_u(Q))(R_u(P)\cap M)(R_u(P)\cap R_u(Q))\] is a 
normal subgroup of $P \cap Q$
(note that $[L\cap R_u(Q),R_u(P)\cap M]\subseteq R_u(P)\cap R_u(Q)$), 
and $V$, being constructible, is closed.
Now $R_u(P)\cap R_u(Q)$ is unipotent and 
$V/(R_u(P)\cap R_u(Q))\cong (L\cap R_u(Q))\times (R_u(P)\cap M)$ 
is unipotent, so $V$ is unipotent, and $V\subseteq G^0$, 
as $R_u(P)$ and $R_u(Q)$ are contained in $G^0$.  
As $V$ is normalized by $T$, it is 
connected thanks to \cite[Prop.\ 14.4(2a)]{borel}.  
Thus $V\subseteq R_u(P\cap Q)$ and 
$R_u(P\cap Q)= ((L\cap M)\cap R_u(P\cap Q))V$.
Now $L\cap M = L_\lambda(L_\mu)$ is an R-Levi subgroup of $G$,
by part (i),  
so $(L\cap M)^0 = (L\cap M) \cap G^0$ is a connected reductive group and 
$(L\cap M)\cap R_u(P \cap Q)= \{ 1 \}$.
We deduce that $V = R_u(P\cap Q)$, as required.
\end{proof}

The next result follows immediately from part (ii) of Lemma \ref{lem:nonconn}.

\begin{cor}
\label{cor:minpar}
 Let $P$ be an R-parabolic subgroup of $G$ with R-Levi subgroup $L$, and let 
$H$ be a closed subgroup of $L$.  Then $P$ is minimal amongst the 
R-parabolic subgroups of $G$ 
containing $H$ if and only if $H$ is $L$-irreducible.
\end{cor}

\begin{cor}
\label{cor:alltor}
 Let $P$ be an R-parabolic subgroup of $G$ and let $T$ be a 
maximal torus of $P$.  Then there exists $\lambda\in Y(T)$ 
such that $P = P_\lambda$.  Moreover, $L_\lambda$ is the unique 
R-Levi subgroup of $P$ that contains $T$. 
\end{cor}

\begin{proof}
Choose $\lambda\in Y(G)$ such that $P=P_\lambda$.  Since maximal 
tori in $P$ are $P$-conjugate, there exists $x\in P$ such that 
$x\lambda(k^*)x^{-1}\subseteq T$.  
We have $T\subseteq L_{x{\cdot}\lambda}$, 
and it is easily checked that $P_{x{\cdot}\lambda} = x P_\lambda x^{-1}= P_\lambda=P$.  
If $L$ is another R-Levi subgroup of $P$ containing $T$, 
then it follows from Lemma \ref{lem:nonconn}(iii), 
setting $P=Q$ and $M=L_\lambda$, that $L=L_\lambda$.
\end{proof}

\begin{cor}
\label{cor:leviinc}
Let $P$ and $Q$ be R-parabolic subgroups of $G$ with $P\subseteq Q$ 
and let $L$ be an R-Levi subgroup of $P$.  Then there is a unique 
R-Levi subgroup $M$ of $Q$ such that $L\subseteq M$.
\end{cor}

\begin{proof}
Choose $\mu\in Y(G)$ such that $P=P_\mu$ and $L=L_\mu$.  
Choose a maximal torus $T$ of $Q$ with $\mu(k^*)\subseteq T$.  
By Corollary \ref{cor:alltor}, we can find $\lambda\in Y(T)$ such that 
$Q = P_\lambda$.  Applying Lemma \ref{lem:nonconn}(iii) 
to the R-parabolic subgroups $P$ and $Q$ and their respective R-Levi 
subgroups $L_\mu$ and $L_\lambda$, yields $L_\mu\subseteq L_\lambda$, 
so we can take $M = L_\lambda$.  Uniqueness follows from 
Corollary \ref{cor:alltor}, as any R-Levi subgroup contains 
a maximal torus of $G$.
\end{proof}

\begin{cor}
\label{cor:leviconj}
 Let $P$ be an R-parabolic subgroup of $G$, 
let $L$ be an R-Levi subgroup of $P$, 
and let $M$ be a closed subgroup of $P$.  
Then $M$ is an R-Levi subgroup of $P$ 
if and only if $M$ is $R_u(P)$-conjugate to $L$.
\end{cor}

\begin{proof}
Write $P = P_\lambda$, $L = L_\lambda$ for some $\lambda\in Y(G)$.  
Clearly, if $u\in R_u(P)$, then $uL_\lambda u^{-1}=L_{u{\cdot}\lambda}$ 
and $P_{u{\cdot}\lambda}= uP_\lambda u^{-1}= P_\lambda$, 
so $uLu^{-1}$ is an R-Levi subgroup of $P$.  
Conversely, suppose that $M$ is an R-Levi subgroup of $P$, say $M = L_\mu$ 
with $P=P_\mu$ for some $\mu\in Y(G)$.  
Since maximal tori of $P$ are $P$-conjugate and $P=R_u(P)M$, 
the R-Levi subgroups $L$ and $uMu^{-1}$ contain a common maximal torus 
for some $u\in R_u(P)$.  
Now Corollary \ref{cor:alltor} implies that $uMu^{-1}=L$.
\end{proof}

\begin{cor}
\label{cor:leviconn}
 Let $P$ be an R-parabolic subgroup of $G$ and let $L$ be a 
Levi subgroup of $P^0$.  Then there exists $\lambda\in Y(G)$ 
such that $P=P_\lambda$ and $L=L_\lambda^0$.
\end{cor}

\begin{proof}
Let $T$ be a maximal torus of $L$.
By Corollary \ref{cor:alltor}, there exists $\lambda \in Y(T)$
such that $P = P_\lambda$.
As $L$ and $L_\lambda^0$ are both Levi subgroups of $P^0$
containing the maximal torus $T$, Corollary \ref{cor:alltor}
implies that $L = L_\lambda^0$.
\end{proof}

The following result is the 
generalization of \cite[Prop 4.4(b)]{boreltits} for non-connected $G$.

\begin{cor}
\label{cor:parintersect}
If $P$ and $Q$ are R-parabolic subgroups of $G$, 
then $(P\cap Q)R_u(Q)$ is an R-parabolic subgroup of $G$.  
Moreover, $(P\cap Q)R_u(Q)=Q$ if and only if 
$P$ contains an R-Levi subgroup of $Q$.
\end{cor}

\begin{proof}
By a standard result for connected groups, \cite[2.4]{boreltits}, $P\cap Q$ 
contains a maximal torus $T$ of $G$.  
We can write $P=P_\mu$ and $Q = P_\lambda$ 
for some $\lambda,\mu\in Y(T)$, by Corollary \ref{cor:alltor}.  
We have
$(P\cap Q)R_u(Q) = (P_\mu\cap L_\lambda)\ltimes R_u(P_\lambda)
= P_\mu(L_\lambda)\ltimes R_u(P_\lambda)$, 
and this is an R-parabolic subgroup of $G$ by Lemma \ref{lem:nonconn}(ii).  
From the first equality we see that if $(P\cap Q)R_u(Q)=Q$, 
then $L_\lambda\subseteq P$.  
Conversely, it is clear that if $P$ contains an R-Levi subgroup of $Q$, 
then $(P\cap Q)R_u(Q)=Q$.
\end{proof} 

\begin{cor}
If $S$ is a torus of $G$, then $C_G(S)$ is an R-Levi subgroup of $G$.
\end{cor}

\begin{proof}
If $S$ is central in $G$, then $C_G(S)=G=P_\lambda$, where $\lambda$ 
is the trivial cocharacter of $G$.  
Otherwise we can find $\lambda\in Y(S)$ such that 
$\lambda(k^*)\not\subseteq Z(G)$.  Then $L_\lambda$ is a proper 
non-connected reductive subgroup of $G$
containing $S$, and $C_G(S)=C_{L_\lambda}(S)$.  
By noetherian induction on closed subgroups of $G$, we can assume that $C_G(S)$ 
is an R-Levi subgroup of $L_\lambda$, say $C_G(S)=L_\mu(L_\lambda)$ 
for some $\mu\in Y(L_\lambda)$.  We have $L_\mu(L_\lambda)=L_{m\lambda+\mu}$ 
for some $m\in {\mathbb N}$ by Lemma \ref{lem:nonconn}(i), as required.
\end{proof}

\begin{lem}
\label{lem:opposite}
Let $P$ be an R-parabolic subgroup of $G$ with R-Levi subgroup $L$.  
Then there exists a unique R-parabolic subgroup $P^-$ of $G$ such 
that $P\cap P^- = L$.
\end{lem}

\begin{proof}
Write $P=P_\lambda$, $L=L_\lambda$ for some $\lambda\in Y(G)$.  
It is simple to check that $P_\lambda\cap P_{-\lambda}=L_\lambda$.  
Uniqueness follows easily from uniqueness in the connected case.
\end{proof}

\begin{lem}
\label{lem:GvsG0}
\begin{itemize}
\item[(i)] 
Let $H$ be a closed subgroup of $G^0$.  
Then $H$ is $G$-completely reducible if and only if 
$H$ is $G^0$-completely reducible.
\item[(ii)]  
If $H$ is a $G$-completely reducible 
subgroup of $G$, then $H\cap G^0$ is $G^0$-completely reducible.
\end{itemize}
\end{lem}

\begin{proof}
(i). Suppose that $H$ is $G$-cr.  Let $P$ be a parabolic 
subgroup of $G^0$ with $H\subseteq P$.  Then $H$ is contained 
in $N_G(P)$, which is an R-parabolic subgroup of $G$ by 
Proposition \ref{prop:normpar}.  Since $H$ is $G$-cr, 
there is an R-Levi subgroup $L$ of $N_G(P)$ with $H\subseteq L$, 
so $H\subseteq L^0$, which is a Levi subgroup of $N_G(P)^0=P$.  
Conversely, suppose that $H$ is $G^0$-cr.  Let $P$ be an R-parabolic 
subgroup of $G$ with $H\subseteq P$.  Then $H\subseteq P^0$, so there 
exists a Levi subgroup $L$ of $P^0$ with $H\subseteq L$, as $H$ is $G^0$-cr.  
By Corollary \ref{cor:leviconn}, we have $L=L_\lambda^0$ 
for some $\lambda\in Y(G)$ with $P=P_\lambda$, so we 
are done.

(ii). If $H$ is $G$-cr, then $H\cap G^0$ is $G$-cr, by the non-connected 
version of Theorem \ref{Serrequestion}.  Now apply part (i).
\end{proof}

\begin{exmp}
\label{ex:finitenoncr}
In general, the converse of Lemma \ref{lem:GvsG0}(ii) is false.  
For example, let $C_p = \langle a\,|\, a^p\rangle$ 
be the cyclic group of order $p = \Char k$, and let $G = C_p \times \SL_2(k)$.
Let $\gamma= \twobytwo{1}{1}{0}{1}$, and set $H=\langle a\gamma\rangle$.  
Then $H\cap G^0 = \{1\}$ is $G^0$-cr and $H\subseteq P_\lambda$, 
where $\lambda\in Y(G)$ is given by 
$\lambda(t)=  \twobytwo{t}{0}{0}{t^{-1}}$, 
but $H$ is not contained in any R-Levi subgroup of $P_\lambda$.
\end{exmp}

\begin{lem}
\label{lem:normalcomplement}
Let $N$ be a closed normal subgroup of $G$.  
Then there exists a closed subgroup $M$ of $G$ 
such that $G = MN$, $M\cap N$ is finite, and $M^0$ commutes with $N^0$.
\end{lem}

\begin{proof}
Clearly, we can assume that $N$ is connected.  
Thus $[N,N]$ is the product of certain simple factors of $G^0$.  
Let $M_1$ be the product of the remaining simple factors of $G^0$.  
As $G$ permutes the simple factors of $G^0$ and $[N,N]$ 
is normalized by $G^0$, $M_1$ is also normalized by $G$; 
moreover, $N$ and $M_1$ commute.
By the proof of \cite[Prop.\ 3.2]{martin1}, there exists 
a finite subgroup $F$ of $G$ such that $G=F G^0$.  
The torus $S:=(Z(G^0)\cap N)^0$ is normal in $G$, so Lemma 2.1 
of \cite{martin1} implies that there exists a subtorus $S_1$ of 
$Z(G^0)^0$ such that $S_1$ is normal in $G$,  
$SS_1=Z(G^0)^0$, and $S\cap S_1$ is finite.  
Set $M=F M_1 S_1$.  It is straightforward to check that $M$ 
has the required properties.
\end{proof}

We now indicate which results from the earlier sections hold for 
non-connected groups.  Subsections \ref{sub:git-gcr}--\ref{sub:rationality} 
below deal with the material in Sections \ref{s:prelim}--\ref{s:rationality} 
respectively.  Our convention, except for in Subsection \ref{sub:building}, 
is that all of the results 
and discussion go through to the non-connected case unless otherwise 
stated below, apart from obvious exceptions such as examples involving 
connected groups.  If the proof of a result is significantly different 
in the non-connected case, then we describe the necessary modifications.

\subsection{Geometric invariant theory and basic properties of $G$-cr subgroups}
\label{sub:git-gcr}

(See Section \ref{s:prelim}).  
Everything in Subsection \ref{sub:git}
goes over to the non-connected case 
with minor modifications, 
cf.\ \cite[\S 4]{martin1}, \cite[\S 1]{martin2}; 
for example, if $G$ acts on an affine variety $X$ and $x\in X$, 
then $G\cdot x$ is closed if and only if $G^0\cdot x$ is closed, 
so the Hilbert--Mumford Theorem \ref{thm:HMT} 
generalizes immediately to non-connected $G$.  
A proof that a strongly reductive subgroup of 
$G$ is non-connected reductive may be found in \cite[\S 6]{martin1}.  
It follows from this and the non-connected version of 
Theorem \ref{mainthm} that a $G$-cr subgroup of $G$ 
is non-connected reductive.  The non-connected version of 
Lemma \ref{linred-cr} follows from the proof 
of \cite[Prop.\ 6.6]{martin1}.  
In Lemma \ref{lem:fgapprox}, the finiteness of 
the number of conjugacy classes of R-parabolic 
subgroups and R-Levi subgroups follows from 
\cite[Prop.\ 5.2(e)]{martin1} and Corollary \ref{cor:leviconj}.
  
Lemma \ref{lem:isogenynew} holds in 
the non-connected case.  
To see this, note that for $f : G_1 \to G_2$ an isogeny
of non-connected reductive groups, 
$\ker{f}$ is a finite normal subgroup of $G_1$, 
so it is centralized by $G_1^0$; this implies that 
$\ker{f}$ is contained in any R-Levi subgroup of 
$G_1$, so part (ii) follows from part (i) as in the connected case.  
The argument given shows that $f(L_\lambda)=L_\mu$ and that $L_\lambda$ 
is proper if and only $L_\mu$ is proper.  
We have $f(P_\lambda^0)=P_\mu^0$ from the connected case, 
so $f(P_\lambda)=f(L_\lambda P_\lambda^0)= L_\mu P_\mu^0= P_\mu$, 
and (i) holds.  Part (iv) follows from the connected case, because 
$R_u(P_\lambda)$ and $R_u(P_\mu)$ are connected.

Lemma \ref{lem:epimorphisms} is more complicated 
for non-connected $G$: 
the underlying problem is that a normal torus of $G$ need not be central.  
Part (i) of Lemma \ref{lem:epimorphisms} holds in the 
non-connected setting, but part (ii)(b) does not for 
the $G_2$-ir and $G_2$-ind cases.  
For let $G=C_2 \ltimes k^*$, where $C_2=\langle a\,|\,a^2\rangle$ 
and $a$ acts on $k^*$ by $a \cdot t = t^{-1}$, 
and let $f : G\ra G/G^0$ be the canonical projection.  
Then clearly $G^0$ is contained in a proper R-Levi subgroup of $G$, 
but $f(G^0) = \{1\}$ is trivially $G/G^0$-ir.

The rest of Lemma \ref{lem:epimorphisms} holds 
in the non-connected case, but we need a different proof.  
Let $f:G_1\ra G_2$ be a surjective homomorphism of 
non-connected reductive groups, with kernel $N$.  
Let $M\subseteq G_1$ satisfy the conclusion of 
Lemma \ref{lem:normalcomplement} with respect to $N$.  
Let $g:M\ra G_2$ be the restriction of $f$; note that $g$ is an isogeny.

First we need a generalization of Lemma \ref{lem:isogenynew}.

\begin{lem}
\label{lem:epiimage}
 Let $\lambda\in Y(M)$ and set $\mu=f\circ \lambda$.  Then
 \begin{itemize}
  \item[(i)] $f(P_\lambda(G_1))=P_\mu(G_2)$,\ 
$f(L_\lambda(G_1))=L_\mu(G_2)$;
  \item[(ii)] $f^{-1}(P_\mu(G_2))=P_\lambda(G_1)$,\ 
$f^{-1}(L_\mu(G_2))=L_\lambda(G_1)$;
  \item[(iii)] $R_u(P_\lambda(G_1))= R_u(P_\lambda(M))$;
  \item[(iv)] $f(R_u(P_\lambda(G_1)))= R_u(P_\mu(G_2))$.
 \end{itemize}
\end{lem}

\begin{proof}
By the non-connected version of Lemma \ref{lem:isogenynew}, 
we can assume that $N$ is connected.  Then $N$ commutes with 
$M^0$, so $N$ lies in $L_\lambda(G_1)$.  
To prove (ii), therefore, it suffices to prove (i).  Clearly, 
$P_\lambda(G_1)= P_\lambda(M)N$, so 
$f(P_\lambda(G_1))= f(P_\lambda(M))= g(P_\lambda(M))= P_\mu(G_2)$, 
by the non-connected version of Lemma \ref{lem:isogenynew}.  
A similar argument gives $f(L_\lambda(G_1))=L_\mu(G_2)$, as required.
If $u\in R_u(P_\lambda(G_1))$, then, writing $u = mn$ with 
$m \in M$, $n \in N$, we have $1 = c_\lambda(u) = c_\lambda(m)c_\lambda(n) =
c_\lambda(m) n$, whence $n = c_\lambda(m)\inverse \in M$. 
Thus $u \in  R_u(P_\lambda(G_1)) \cap M =  R_u(P_\lambda(M))$,
and part (iii) is proved.
Part (iv) now follows from (iii) and 
the non-connected analogue of Lemma \ref{lem:isogenynew}(iv).
\end{proof}

Now suppose that $H_1$ is $G_1$-cr.  
Suppose that $f(H_1)\subseteq P_\mu(G_2)$, 
where $\mu\in Y(G_2)$.  Write $\mu=f\circ \lambda$, 
where $\lambda\in Y(M)$.  
Then $f^{-1}(f(H_1))\subseteq f^{-1}(P_\mu(G_2))= P_\lambda(G_1)$, 
by Lemma \ref{lem:epiimage}(ii), so $H_1\subseteq P_\lambda(G_1)$.  
Since $H_1$ is $G_1$-cr, $H_1$ lies in an R-Levi subgroup of 
$P_\lambda(G_1)$, so by Corollary \ref{cor:leviconj}, we have 
$H_1\subseteq L_{u\cdot \lambda}(G_1)$ for some $u\in R_u(P_\lambda(G_1))$.  
Now $u\in R_u(P_\lambda(M))$, by Lemma \ref{lem:epiimage}(iii), 
so $u\cdot \lambda\in Y(M)$.  
Thus $f(H_1)\subseteq f(L_{u\cdot \lambda}(M))= L_{f(u)\cdot \mu}(G_2)$, 
by Lemma \ref{lem:epiimage}(i); 
moreover, $P_{f(u)\cdot \mu}(G_2) = 
f(P_{u\cdot \lambda}(G_1))= f(P_\lambda(G_1))= P_\mu(G_2)$ 
(Lemma \ref{lem:epiimage}(i)), 
so $L_{f(u)\cdot \mu}(G_2)$ is 
an R-Levi subgroup of $P_\mu(G_2)$.  It follows that $f(H_1)$ is $G_2$-cr.

Suppose that $H_1$ is $G_1$-ind.  If $f(H_1)\subseteq L_\mu(G_2)$ 
for some $\mu\in Y(G_2)$ with $L_\mu(G_2)\neq G_2$ then, picking 
$\lambda\in Y(M)$ such that $\mu=f\circ \lambda$, we have 
$H_1\subseteq f^{-1}(f(H_1))\subseteq f^{-1}(L_\mu(G_2))= 
L_\lambda(G_1)$, by Lemma \ref{lem:epiimage}(ii).  
By Lemma \ref{lem:isogenynew}(iii), $L_\lambda(M)\neq M$, 
so $\lambda(k^*)$ does not centralize $M$, so $\lambda(k^*)$ 
does not centralize $G_1$, so $L_\lambda(G_1)\neq G_1$.  
But this contradicts the $G_1$-indecomposability of $H_1$.  
We deduce that $f(H_1)$ is $G_2$-ind.

It now follows that if $H_1$ is $G_1$-ir then $f(H_1)$ is $G_2$-ir.
This completes the proof of Lemma \ref{lem:epimorphisms}(ii)(a) in 
the non-connected case.

Now we show that if $f$ is non-degenerate and $H_1$ is a closed 
subgroup of $G_1$ such that $f(H_1)$ is $G_2$-cr, then $H_1$ is $G_1$-cr.  
By Lemma \ref{lem:isogenynew}(i), we can assume that $N$ 
is a torus.  Thus $N$ lies in every R-Levi subgroup of $G_1$, so 
without loss we can assume that $H_1\subseteq M$.  
Let $\lambda\in Y(G_1)$ with $H_1\subseteq P_\lambda(G_1)$.  
Then $f(H_1)\subseteq P_\mu(G_2)$, where $\mu:= f\circ \lambda$.  
As $f(H_1)$ is $G_2$-cr, there is an R-Levi subgroup $L$ of $P_\mu(G_2)$ 
with $f(H_1)\subseteq L$.  

Choose $\sigma\in Y(M)$ such that $g\circ \sigma=\mu$.  
We have $\lambda(k^*)\subseteq \sigma(k^*)N$; 
since $N\subseteq Z(G_1^0)$, there exists $\tau\in Y(N)$ 
such that $\lambda= \sigma+\tau$.  
By Lemma \ref{lem:epiimage}(iv) and Corollary \ref{cor:leviconj}, 
there exists $u\in R_u(P_\sigma(M)) = R_u(P_\lambda(M)) \subseteq 
R_u(P_\lambda(G_1))$ such that 
$f(H_1)\subseteq L_{f(u){\cdot}\mu}(G_2)$.  
Thus, replacing $\lambda$ by $u{\cdot}\lambda$ if necessary, 
we can assume that $f(H_1)\subseteq L_\mu(G_2)$.
We have $H_1\subseteq g^{-1}(L_\mu(G_2))=L_\sigma(M)$, 
by Lemma \ref{lem:isogenynew}(ii).
We have $\lambda\in Y(L_\sigma(G_1))$ and 
\[
H_1\subseteq P_\lambda(G_1)\cap L_\sigma(M) 
\subseteq P_\lambda(G_1)\cap L_\sigma(G_1) 
= P_\lambda(L_\sigma(G_1)) = P_\tau(L_\sigma(G_1)).
\]  
Now $L_\sigma(G_1)^0\subseteq P_\tau(L_\sigma(G_1))$, 
thus $R_u(P_\tau(L_\sigma(G_1)))= \{1 \}$,
and so $P_\tau(L_\sigma(G_1))=L_\tau(L_\sigma(G_1))$.  
Thus $H_1\subseteq L_\tau(L_\sigma(G_1))\subseteq L_\lambda(G_1)$.  
We conclude that $H_1$ is $G_1$-cr.  
This completes the proof of the remaining parts 
of Lemma \ref{lem:epimorphisms} 
in the non-connected case.

\bigskip

For Propositions \ref{rich-stable} and \ref{Rich-closedorbit}, 
see \cite[Prop.\ 8.3]{martin1}.  Proposition \ref{Rich-linearlyreductive} 
for non-connected $G$ follows from the connected case,
Lemma \ref{lem:GvsG0} and the non-connected version of Theorem~\ref{mainthm}.

\subsection{$G$-complete reducibility}
\label{sub:cr}

(See Section \ref{s:cr}).
The proofs of Theorem \ref{mainthm} and Corollary \ref{mainthm-cor} 
go through from the connected case to the non-connected case, 
replacing various results concerning parabolic subgroups with 
their counterparts from Subsection \ref{sub:prelim} 
for R-parabolic subgroups.  For example, we use Lemma \ref{lem:nonconn}(ii) 
instead of \cite[Prop.\ 4.4(c)]{boreltits}, Corollary \ref{cor:parintersect} 
instead of \cite[Prop.\ 4.4(b)]{boreltits} and 
\eqref{eq:intersection} 
instead of \eqref{fac}.

Proposition \ref{ben3} generalizes (with the same proof)
to give the result for non-connected $G$ and $H$.
For the proof of Proposition \ref{redcent}, note that if $M$ is 
a closed subgroup of $G$, then $G/M$ is affine if and only if 
$G^0/(G^0\cap M)$ is affine (cf.\ \cite[\S 1]{rich}).
By virtue of these extensions, we now obtain
Theorem \ref{converse} immediately as a special case
of Proposition \ref{ben3} --- the separate proof in 
Section \ref{s:cr} is not needed.

The first assertion of Corollary \ref{kempf1} fails in the 
non-connected case, by the counterexamples below to 
Lemma \ref{lem:irrcentlinred}.  The second assertion, 
however, still holds.  For suppose that $H$ is not $G$-cr.  
By Remark \ref{rem:topfg} and Corollary \ref{cor:closedorbit}, 
we can assume that $H$ is topologically finitely generated, say 
by $h_1,\ldots, h_n$, and that the orbit $G\cdot (h_1,\ldots, h_n)$ 
is not closed in $G^n$.  
By the Hilbert--Mumford--Kempf Theorem (see \cite[Cor.\ 3.5]{kempf}), 
there exists $\lambda\in Y(G)$ such that 
$H\subseteq P_\lambda$, $C_G(H)\subseteq N_G(P_\lambda)$, 
$G\cdot (c_\lambda(h_1), \ldots, c_\lambda(h_n))$ is closed but 
$G\cdot (h_1,\ldots, h_n)$ is not.  
This implies that $P_\lambda\neq L_\lambda$, 
whence $R_u(P_\lambda) \neq \{ 1 \}$,
whence $P_\lambda^0\neq G^0$, whence $N_G(P_\lambda^0)$ 
is a proper R-parabolic subgroup of $G$, by Proposition \ref{prop:normpar}. 
But this is impossible, since $C_G(H)$ is $G$-ir.

An immediate corollary of Proposition \ref{regular-is-cr}
is the following: if $H$ is a subgroup of $G$ containing $G^0$, then $H$ 
is $G$-cr.

Corollary \ref{cor:levi} holds in the non-connected case, 
although the proof given in Section \ref{s:cr} does not work,
because $C_G(S)$ need not be connected.  
Here is an alternative proof.  
Let $L$ be an R-Levi subgroup of $G$.
Recall that $L = C_G(Z(L)^0)$.
Let $K$ be a closed subgroup of $L$ and 
let $S$ be a maximal torus of 
$C_G(K)$ with $Z(L)^0\subseteq S$.  Then $S\subseteq L$, so 
$S$ is a maximal torus of $C_L(K)$ and $C_G(S)= C_L(S)$.  
The desired result now follows from 
the non-connected analogue of Corollary \ref{mainthm-cor}.

We do not know whether Theorem \ref{gcr-regular} 
and Remark \ref{rem:seprem} hold in the non-connected case; the proof of 
\cite[Thm.\ 4.5]{martin3} does not generalize.

We extend the definition of a reductive pair to 
non-connected $G$ and $H$ in the obvious way.  
Example \ref{exmp:redpair} is valid in the non-connected case: 
for $H$ is generated by $H^0$ and $N_H(T)$, and $N_H(T)$ 
permutes $\Psi(G,T)\setminus\Psi(H,T)$.  
Also Example \ref{ex:adjoint} is valid for non-connected 
$H$ with $H^0$ a simple group of adjoint type.

Lemma \ref{lem:irrcentlinred} fails in the non-connected case: 
e.g., if $G$ is an abelian $p$-group, then $Z(G)$ itself is not 
linearly reductive.  It is not even true in general that $C_G(H)/Z(G)$ 
is linearly reductive: take $G$ to be $C_2\ltimes (M\times M)$, 
where $p=2$, $M$ is a connected simple group and $C_2$ acts 
by permuting the factors, and take $H$ to be the copy of $M$ 
diagonally embedded in $G^0$.

We consider the material in Subsection \ref{sub:ss} for 
connected $G$ only.  Theorem \ref{thm:adjcr}, however, 
holds in the non-connected case.

\subsection{Complete reducibility in $G$ and the building of $G^0$}
\label{sub:building}

(See Section \ref{s:building}).
We will not construct a simplicial complex using the 
R-parabolic subgroups of $G$ for non-connected $G$.  
Instead we regard $G$ as a subgroup of $\Aut X$ ($X:=X(G^0)$):  
we consider the action of $G$ on the spherical building
$X$ induced by the conjugation action of $G$ on $G^0$.
For a subgroup $H$ of $G$, it need no longer be the case that $X^H$
is a subcomplex of $X$. 
However, we can view the subspace $X^H$ as a subcomplex of
the barycentric subdivision of $X$, see \cite[2.3.1]{serre2}. 

It is still possible to ask whether the subspace $X^H$ is 
$X$-completely reducible or not.
The following proposition allows us to extend the building-theoretic
interpretation of $G$-complete reducibility to this situation.

\begin{prop}
\label{nonconnectedbuilding}
Let $H$ be a subgroup of $G$. 
Then $H$ is $G$-completely reducible if and only if every 
simplex in $X^H$ has an
opposite in $X^H$.
\end{prop} 

\begin{proof}
Recall that for any parabolic subgroup $Q$ of $G^0$,
$N_G(Q)^0 = N_{G^0}(Q) = Q$, see the note following
Proposition \ref{prop:normpar}.
Now suppose that $H$ is $G$-cr.
Let $Q$ be a parabolic subgroup of $G^0$ fixed by $H$.
Then $H \subseteq N_G(Q)$, which is an R-parabolic subgroup of $G$,
by Proposition \ref{prop:normpar}.
Since $H$ is $G$-cr, there exists $\lambda\in Y(G)$ with 
$N_G(Q) = P_\lambda$ and 
$H\subseteq L_\lambda = P_\lambda \cap P_{-\lambda}$.
But then $H$ fixes $P_{-\lambda}^0$, which is a
parabolic subgroup of $G^0$ opposite to $P_\lambda^0 = Q$.
Thus every simplex in $X^H$
has an opposite in $X^H$.

Conversely, suppose that every simplex 
in $X^H$ has an opposite in $X^H$.
Let $P$ be an R-parabolic subgroup of $G$ containing $H$.
Then $Q = P^0$ is an $H$-fixed simplex of $X$, so there
exists an $H$-fixed simplex $Q^-$ of $X$ opposite $Q$.
Let $M = Q\cap Q^-$ be the common Levi subgroup of $Q$ and
$Q^-$.
Set $P' = N_G(Q^-)$.  
By Corollary \ref{cor:leviconn}, there exist R-Levi subgroups $L$ and $L'$ 
of $P$ and $P'$ respectively such that $L\cap G^0 = L'\cap G^0 = M$.
Then 
$$
P\cap P' = (L\cap L')(L\cap R_u(P'))(L'\cap R_u(P))(R_u(P)\cap R_u(P')),
$$
by Lemma \ref{lem:nonconn}(iii).
However, 
$L\cap R_u(P') = M\cap R_u(Q^-) = \{1\}$,
$L'\cap R_u(P) = M\cap R_u(Q) = \{1\}$, and 
$R_u(P) \cap R_u(P') =  R_u(Q) \cap R_u(Q^-) = \{1\}$.
Thus $H \subseteq P\cap P' = L\cap L' \subseteq L$,
so $H$ is $G$-cr, as claimed.
\end{proof}

We now see that part of Theorem \ref{gcr-building} holds 
for non-connected $G$: a subgroup $H$ of $G$ is $G$-cr if and only if the
fixed point subspace $X^H$ is $X$-cr.
We are grateful to J.-P. Serre for drawing our attention
to this extension.

Lemma \ref{parabolics-cent} holds for non-connected
$G$ and also holds with $C_G(H)^0$ replaced everywhere
by $C_G(H)$.
Moreover, by the arguments in the proof of 
Proposition \ref{fixbuilding-regular}, the following is also true:
if $H$ is any closed (not necessarily connected) $G$-cr 
subgroup of $G$, 
then R-parabolic subgroups of $HC_G(H)^0$ containing $H$ 
correspond bijectively to R-parabolic subgroups of $C_G(H)^0$, 
and the same is true replacing $C_G(H)^0$ everywhere with $C_G(H)$.

When the action of $H$ is
completely reducible and type-preserving,
then  $X^H$ is a subcomplex of $X$ and
the results due to M\"uhlherr also hold: i.e., $X^H$
is a building if and only if it is thick, and
the thickening of $X^H$ is isomorphic to 
$X(C_{G^0}(H)^0) = X(C_G(H)^0)$.

\subsection{Rationality questions}
\label{sub:rationality}

(See Section \ref{s:rationality}).  
Theorem \ref{thm:perfectcr} holds in the non-connected case, 
by \cite[Prop.\ 1.2]{borel}: if $M$ is a linear algebraic 
group defined over $k'$, then $M^0$ is also defined over $k'$.
The existence of a $\Gal(\overline{k}/k)$-invariant length function
on $Y(G)$ follows from the arguments of \cite[\S 4]{kempf}
and \cite[\S 1]{martin2}.
We use Corollary \ref{cor:alltor} instead of \cite[Cor.\ 8.4.4]{spr2}.  
Note that \cite[Thm.\ 4.2]{kempf} holds for non-connected $G$ 
--- cf.\ the discussion of the Hilbert--Mumford 
Theorem in Subsection \ref{sub:git-gcr} --- 
and the argument of \cite[Prop.\ 2.2]{LMS} 
also works in the non-connected case.


\bigskip {\bf Acknowledgements}:
The first author acknowledges funding by the EPSRC. 
We would like to thank  M.W.\ Liebeck, 
G.\  McNinch, B.\  M\"uhlherr, and G.M.\ Seitz
for helpful discussions.
We are especially grateful to J.-P.\  Serre for his 
many helpful comments on earlier versions of the paper and
in particular for pointing out the connection between 
Theorem \ref{mainthm} and his Theorem \ref{gcr-building}.

\smallskip


\end{document}